\documentclass[12pt]{amsart}
\usepackage{amssymb}
\usepackage[dvips]{graphics}

\ifx\emph\undefined\let\emph\it\fi
\ifx\mathrm\undefined\let\mathrm\rm\fi
\ifx\mathbf\undefined\let\mathbf\bf\fi
\ifx\mathfrak\undefined\let\mathfrak\frak\else\let\frak\mathfrak\fi
\ifx\mathcal\undefined\let\mathcal\cal\else\let\cal\mathcal\fi
\ifx\mathbb\undefined\let\mathbb\Bbb\else\let\Bbb\mathbb\fi

\newcommand{\g}{{{\mathfrak g}\,}}

\newcommand{\n}{{{\mathfrak n}}}
\newcommand{\h}{{{\mathfrak h\,}}}

\newcommand{\Z}{{\mathbb Z}}

\newcommand{\C}{{\mathbb C}}

\newcommand{\Ref}[1]{{(\ref{#1})}}

\newcommand{\la}{\lambda}

\newcommand{\dontprint}[1]

\newcommand{\nc}{\newcommand}

\newcommand{\bs}{\boldsymbol}
\nc{\Wr}{{ {\rm Wr}}}
\newcommand{\beq}{\begin{equation}}
\newcommand{\eeq}{\end{equation}}
\newcommand{\sing}{{\rm Sing}\,}

\newcommand{\bean}{\begin{eqnarray}}
\newcommand{\eean}{\end{eqnarray}}
\newcommand{\be}{\begin{displaymath}}
\newcommand{\ee}{\end{displaymath}}
\newcommand{\bea}{\begin{eqnarray*}}
\newcommand{\eea}{\end{eqnarray*}}

\def\Deg{\deg}

\def\dpar{\partial}

\def\Ml{M_{\bs \La}[\bs l]}

\def\slt{\mathfrak{sl}_2}
\def\slth{\mathfrak{sl}_3}
\def\slg{\mathfrak{sl}}
\def\glg{\mathfrak{gl}}
\def\p{\partial}

{\relax}

\newtheorem%
{thm}{Theorem}[section]
\newtheorem%
{proposition}[thm]{Proposition}
\newtheorem%
{lemma}[thm]{Lemma}
\newtheorem%
{lemmadef}[thm]{Lemma-Definition}
\newtheorem%
{corollary}[thm]{Corollary}
\newtheorem%
{conjecture}[thm]{Conjecture}

\nc{\al}{\alpha}
\nc{\om}{\omega}
\nc{\La}{\Lambda}
\nc{\un}{U(\n_-)}
\nc{\px}{\frac d {d x}}
\nc{\PCr}{{ \bs P (\C[x])^r }}

\nc{\on}{\operatorname}

\begin{document}

\title[On separation of variables in the $\slth$ Gaudin model]
{On the new form of Bethe ansatz equations and
separation of variables in the $\slth$ Gaudin model}

\author[E.\,Mukhin, V.\,Schechtman, V.\,Tarasov, and A.\,Varchenko]
{E.\,Mukhin$\,^{*,1}$, V.\,Schechtman$\,^{**}$, V.\,Tarasov$\,^{*,\star,2}$,
\and A.\,Varchenko$\,^{\star\star,3}$}

\thanks{${}^1$\ Supported in part by NSF grant DMS-0601005}
\thanks{${}^2$\ Supported in part by RFFI grant 05-01-00922}
\thanks{${}^3$\ Supported in part by NSF grant DMS-0555327}

\maketitle

\centerline{\it ${}^*$Department of Mathematical Sciences,
Indiana University -- Purdue University,}
\centerline{\it Indianapolis, 402 North Blackford St, Indianapolis,
IN 46202-3216, USA}
\smallskip
\centerline{\it $^\star$St.\,Petersburg Branch of Steklov Mathematical
Institute,}
\centerline{\it Fontanka 27, St.\,Petersburg, 191023, Russia}
\smallskip

\centerline{\it ${}^{**}$
Laboratoire Emile Picard, UFR MIG,
Universit\'e Paul Sabatier,}

\centerline{\it 31062 Toulouse, France}

\centerline{\it ${}^{\star\star}$Department of Mathematics, University of
North Carolina at Chapel Hill,} \centerline{\it Chapel Hill, NC
27599-3250, USA}

\bigskip
\begin{center}
{\it Dedicated to V.I.\,Arnold on the occasion
of his $70^{\text{th}}$ birthday.}
\end{center}
\bigskip

\medskip

\thispagestyle{empty}

\begin{abstract}
A new form of Bethe ansatz equations is introduced. A version
of a separation of variables for the quantum $\slth$ Gaudin model is presented.
\end{abstract}

\section{Introduction}
The separation of variables for the quantum $\slt$ Gaudin model was constructed
by Sklyanin in \cite{Sk}. In this paper, we give an analogue of Sklyanin's
construction for the Lie algebra $\slth$.

We were inspired by Stoyanovsky's paper \cite{St}, in which the author uses
Sklyanin's change of variables to establish a relation between the $\slt$
Knizhnik-Zamolodchikov equations \cite{KZ} and the
Belavin-Polyakov-Zamolodchikov equations \cite{BPZ},
and to construct integral formulae for solutions to the BPZ equations.

\subsection{}
The paper is organized as follows. In Section \ref{sec KZ}, we
introduce notations and recall the definition of the KZ equations and
Gaudin model.

\subsection{}

In Section \ref{MasteR}, we recall the definition of the master
function and canonical weight function. In Theorem \ref{resh}, we
recall the main fact of the Bethe ansatz method: if $\bs t$ is a critical
point of the master function, then the value at $\bs t$ of the
canonical weight function is an eigenvector of the Gaudin
Hamiltonians, see \cite{RV}, cf. \cite{Ba1, Ba2}.

The critical point equations for the master function
are called also the Bethe ansatz equations for the Gaudin model.

\subsection{}

In Section \ref{sec Different}, we discuss different forms of the
Bethe ansatz equations for the Lie algebras $\slt$ and $\slth$.

For $\slt$, there are two ways to describe solutions to the Bethe ansatz
 equations.
The original way: a solution is a collection of numbers
$\bs t = (t^{(1)}_1, \dots,t^{(1)}_{l_1})$,
satisfying the critical point equations \Ref{BAE2}.
The second way:
given the numbers $z_1,\dots,z_n, (\La_i,\al_1), (\La_i,\La_j)$,
a solution is a polynomial
$P(x)=\prod_{i=1}^{l_1}(t^{(1)}_i-x)$ of one variable and numbers
$\mu_1,\dots,\mu_n$, $\mu_1+\dots +\mu_n=0$,
satisfying the differential equation
\bean
\label{eqn2''}
P''\ -\ P' \sum_{i=1}^n\frac{(\La_i,\al_1)}{x-z_i}
\ +\ P \sum_{i=1}^n \frac 1{x-z_i} \,
(\,\mu_i - \sum_{j\neq i}
\frac{ (\La_i,\La_j)}{z_i-z_j}\,)\ = 0\ .
\eean

The equivalence of two descriptions
is a classical fact which goes back to Stieltjes,
see \cite{Sti} and Sec. 6.8 in \cite{Sz}.

For $\slth$, there are four ways to describe solutions
to the Bethe ansatz equations.
In this introduction, we mention only two of the four:
the original way and the new way suggested in this paper.

The original way: a solution is a collection of numbers
$\bs t = (t^{(1)}_1, \dots,t^{(1)}_{l_1},t^{(2)}_1, \dots,t^{(2)}_{l_2})$,
satisfying the critical point equations
\Ref{BAE3}. The new way:
given the numbers $z_1$,
\dots , $z_n$, $(\La_i,\al_1)$, $(\La_i,\al_2)$,
$(\La_i,\La_j)$,
a solution is a pair
$$
P_1(x)\ =\ \prod_{i=1}^{l_1}(t^{(1)}_i-x)\ ,
\qquad
P_2(x)\ =\ \prod_{i=1}^{l_2}(t^{(2)}_i-x)
$$
of polynomials of one variable
and a set of numbers
$\mu_1, \dots , \mu_n$, $\mu_1 + \dots + \mu_n = 0$,
satisfying the differential equation
\bean
\label{BAE new''}
&&
\phantom{aaaaaaa}
P_1''P_2 - P_1'P_2' + P_1P_2'' -
P_1'P_2 \sum_{i=1}^n \frac{(\La_i,\al_1)}{x-z_i}
-
P_1P_2' \sum_{i=1}^n \frac{(\La_i,\al_2)}{x-z_i}
+
\phantom{aaaaaaaaaa}
\\
&&
\phantom{aaaaaaaaaaaaaaaaaaaaaaaaaaaaa}
+ P_1 P_2 \sum_{i=1}^n \frac{1}{x-z_i}
(\,\mu_i - \sum_{j\neq i}
\frac{ (\La_i,\La_j)}{z_i-z_j}\,)\ = 0\ ,
\notag
\eean
see Theorem \ref{new thm}. This theorem
is the first main result of our paper.

\subsection{}

In Section \ref{On separation for slt}, we describe Sklyanin's
separation of variables for the $\slt$ Gaudin model, following the exposition
in \cite{St}.

The general problem in the Gaudin model
is to find eigenfunctions of the
commuting Gaudin Hamiltonians $H_i(\bs z),\, i=1,\dots,n$,\
where
$H_1(\bs z)+\dots+ H_n(\bs z)=0$.
In suitable coordinates,
the Hamiltonians are differential operators acting on
polynomials in $n$ variables
$x^{(1)},\dots,x^{(n)}$.

The eigenfunction equations are
$$
H_i(\bs z) \,F(x^{(1)},\dots,x^{(n)})\ =\ \mu_i\,
F(x^{(1)},\dots,x^{(n)})\ ,
\qquad
i=1,\dots,n\ ,
$$
where the eigenvalues satisfy the equation $\mu_1+\dots +\mu_n=0$.
The famous Sklyanin's change of variables \Ref{change}
from variables $x^{(1)},\dots,x^{(n)}$ to new variables
$u, y^{(1)},\dots, y^{(n-1)}$ transforms the eigenfunction
equations to the following equations for the unknown polynomial
$F(\bs x(u,\bs y))$,
\bean
\label{famous}
\phantom{aaaa}
\bigl(- \p^2_{y^{(j)}}\
+\
\sum_{i=1}^n \, \frac{(\La_i,\al_1)}{y^{(j)}-z_i}\, \p_{y^{(j)}}
\ +\
\sum_{i=1}^n \frac 1{y^{(j)}-z_i} \,
(\,-\mu_i + \sum_{k\neq i}
\frac{(\La_i,\La_k)}{z_i-z_k}\,)\, \Bigr)\, F\ =\ 0\ ,
\eean
$j=1,\dots,n-1$, see \cite{Sk, St} and Theorem \ref{Sklyanin thm}.

This Sklyanin's statement has three interesting features.

The first is that the variables have separated: the $j$-th
equation depends on variable $y^{(j)}$ only and
does not depend on $u$ and other variables $y^{(i)}$.

The second feature is that the differential operator is the same in all equations.
Namely, the differential operator for an index $j'$ can
be obtained from the differential operator for the index $j$ by
replacing $y^{(j)}$ with $y^{(j')}$.

The third feature is that the differential operator in Sklyanin's
equations is the same as the differential operator in the Bethe ansatz equation
\Ref{eqn2''}.

\medskip
Sklyanin and Stoyanovsky also consider the canonical weight
function. In Sklyanin's variables, the canonical weight function takes the form
\bea
\Psi(\bs t, \bs z, u, \bs y)\ =\
u^{l_1}\,
\frac{ P(y^{(1)})\dots P(y^{(n-1)})}
{ P(z_1)\dots P(z_{n})}\ ,
\eea
where
$l_1$ is a nonnegative integer
and
$P(x) = \prod_{i=1}^{l_1}(t^{(1)}_i-x)$.

The canonical weight function, as a function
of $u$, $y^{(1)}$, $\dots,$ $y^{(n-1)}$, is the product
of functions of one variable.
This is another manifestation of separation of variables.

\medskip
Now if one looks for a value of the parameters $\bs t$ such that
$\Psi(\bs t, \bs z, u, \bs y)$,
as a function of $u, \bs y$, becomes an eigenfunction
of the Gaudin operators, then one gets the Bethe ansatz
equation \Ref{eqn2''} for the unknown polynomial $P(x)$.

Hence, if for some $\bs t$, the function
$\Psi(\bs t, \bs z, u, \bs y)$ is an eigenfunction, then $\bs t$
satisfies the Bethe ansatz equations.

\subsection{}

In Section \ref{On separation for slth}, we describe an analog of
Sklyanin's statements for $\slth$.

Again the general problem
is to find eigenfunctions of the
Gaudin Hamiltonians $H_i(\bs z),\ i=1,\dots,n$,\,
where
$H_1(\bs z)+\dots+ H_n(\bs z)=0$.
Now the Hamiltonians are differential operators acting on
polynomials in $3n$ variables
$x^{(1)}_1,x^{(1)}_2, x^{(1)}_3,$
\dots, $x^{(n)}_1,x^{(n)}_2, x^{(n)}_3$.

The eigenfunction equations are
$$
H_i(\bs z)\, F(\bs x)\ =\ \mu_i\, F(\bs x)\ ,
\qquad
i=1,\dots,n\ ,
$$
where the eigenvalues satisfy the equation $\mu_1+\dots +\mu_n=0$.

We make a change
from variables $x^{(1)}_1,x^{(1)}_2, x^{(1)}_3,$
\dots, $x^{(n)}_1,x^{(n)}_2, x^{(n)}_3$ to new variables
$u_1,u_2,u_3$, $y^{(1)}_1,y^{(1)}_2, y^{(1)}_3,$
\dots, $y^{(n-1)}_1,y^{(n-1)}_2, y^{(n-1)}_3$, see formula
\Ref{change slth} analogous to Sklyanin's formula.

We define the degree of a polynomial in $\C[\bs u,\bs y]$
as its degree with respect to variable $u_3$.

On the affine space with coordinates $\bs u,\bs y$, we consider
the affine subspace
$$
\frak D\ =\ \{ (\bs u,\bs y)\ |\ y^{(j)}_1=y^{(j)}_2,\,
j=1,\dots,n-1 \}\ .
$$

It turns out that the Gaudin Hamiltonians have the decomposition
\bea
H_i(\bs z)
= \
\bar H_i(\bs z)_0 + H_i(\bs z)_{>0} + H_i(\bs z)_{<0}\ ,
\eea
where the operator $\bar H_i(\bs z)_0$ preserves the degree,
the operator $ H_i(\bs z)_{>0}$
increases the degree by one, and the operator
$ H_i(\bs z)_{<0}$ decreases the degree by one.
Moreover, it turns out that for any
$F \in \C[\bs u,\bs y]$, the restriction
to $\frak D$ of the polynomial
$H_i(\bs z)_{<0} F$ is zero.

These remarks show that
the eigenfunction equations become
``upper triangular'' with respect to the degree decomposition
after restriction to $\frak D$.

\medskip
Let $F(\bs u ,\bs y)$ be an unknown eigenfunction with eigenvalues
$\mu_1,\dots,\mu_n$. Let
$F(\bs u ,\bs y) = F(\bs u ,\bs y)_0 + F(\bs u ,\bs y)_1 + \dots$
be the degree decomposition. According to our change of variables,
\bea
F(\bs u ,\bs y)_0\ =\ u_1^{l_1} u_2^{l_2}\, f (y^{(1)}_1,y^{(1)}_2,
\dots, y^{(n-1)}_1,y^{(n-1)}_2)\ ,
\eea
where $l_1$, $l_2$ are nonnegative integers and
$f$ is a polynomial depending on variables
$y^{(j)}_1,y^{(j)}_2$,
$j=1,\dots, n-1$, only.

The triangularity remark implies the following
equations for the degree zero part of the eigenfunction,
\bean
\label{degree zero eqnS}
\bigl((\bar H_i(\bs z)_{0} - \mu_i)\, F_0\Bigr) \,
{\Big{\vert}}_{\frak D} \ =\ 0\ , \qquad
i=1,\dots,n\ .
\eean

\subsubsection{}
\label{first thm in intro}
{\bf Theorem.}
{\it After suitable renormalization, equations \Ref{degree zero eqnS}
take the form}
\bean
\label{Main eqn in coordinates}
&&
\Bigl(\,-\, \frac {\p^2f} {\p y^{(j)}_1 \p y^{(j)}_1}
\ +\
\frac {\p^2f} {\p y^{(j)}_1 \p y^{(j)}_2}
\ -\
\frac {\p^2f} {\p y^{(j)}_2 \p y^{(j)}_2}
\ -\
\frac {\p f} {\p y^{(j)}_1}
\sum_{i=1}^n \frac{(\La_i,\al_1)}{y^{(j)}_1 - z_i} \ -
\phantom{\Bigr)}
\\
&&
\phantom{\Bigl( aaaaaaa}
- \
\frac {\p f} {\p y^{(j)}_2}
\sum_{i=1}^n \frac{(\La_i,\al_2)}{y^{(j)}_1-z_i}
\ +\
f \sum_{i=1}^n \frac 1{y^{(j)}_1-z_i} \,
(\,- \mu_i + \sum_{k\neq i}
\frac{(\La_i,\La_k)}{z_i-z_k}\,)
\,\Bigr) {\bigg{\vert}}_{\frak D} \ =\ 0\ .
\notag
\eean
\medskip

This is our second main result, see Theorem \ref{thm on separation}.

As Sklyanin's equations \Ref{famous} for $\slt$, our equations
\Ref{Main eqn in coordinates} have three interesting properties.

The first is that the variables have separated: the $j$-th
equation depends on variables $y^{(j)}_1$, $y^{(j)}_2$ only and
does not depend on $\bs u$ and other variables $y^{(j')}_1$, $y^{(j')}_2$.

The second property is that the differential operator is the same in all equations.
Namely, the differential operator for an index $j'$ can
be obtained from the differential operator for the index $j$ by
replacing $y^{(j)}_1, y^{(j)}_2$ with $y^{(j')}_1, y^{(j')}_2$, respectively.

The third property is that the differential operator in
equations \Ref{Main eqn in coordinates}
is the same as the differential operator in the Bethe ansatz equation
\Ref{BAE new''}.

\medskip

Then we consider the canonical weight
function $\Psi(\bs t, \bs z, \bs u, \bs y)$ and its weight decomposition
$\Psi = \Psi_0 + \Psi_1 + \dots$. It turns our that the degree zero term has the
form
\bean
\label{Canon fn slth}
\Psi(\bs t, \bs z, \bs u, \bs y)_0\ =\
u_1^{l_1} u_2^{l_1}\,
\frac{ \prod_{j=1}^{n-1}\,P_1(y_1^{(j)})
P_2(y_2^{(j)})}
{ \prod_{s=1}^n\, P_1(z_s) P_2(z_s)}\ ,
\eean
where $l_1,l_2$ are nonnegative integers, \
$P_1(x) = \prod_{i=1}^{l_1} (t^{(1)}_i-x)$,\
$P_2(x) = \prod_{i=1}^{l_2} (t^{(2)}_i-x)$.
Moreover, it turns out that the higher degree terms of $\Psi$
are uniquely
determined by the degree zero term $\Psi_0$, see Theorem \ref{slth can fn thm}
and Section \ref{higher terms are determined}.

According to formula \Ref{Canon fn slth},
the degree zero term $\Psi_0$ of the
canonical weight function is the product of functions of one variable.
This is another manifestation of separation of variables.

\medskip

Now assume that one looks for a value of parameters $\bs t$ such that
$\Psi(\bs t, \bs z, \bs u, \bs y)$, as a function of $\bs u, \bs y$,
becomes an eigenfunction of the Gaudin operators. Then, by
Theorem \ref{first thm in intro}, one gets the Bethe ansatz
equation \Ref{BAE new''} for the unknown polynomials $P_1(x)$, $P_2(x)$.

Hence, if for some $\bs t$, the function
 $\Psi(\bs t, \bs z, \bs u, \bs y)$ is an eigenfunction, then $\bs t$
satisfies the Bethe ansatz equations. The precise statement see in
Theorem \ref{eigenvalue slth thm}.

\subsection{} In the appendix we consider the Bethe ansatz equations for the Gaudin model
associated with an arbitrary Kac-Moody algebra of rank $r$. We
introduce a polylinear differential equation for a collection
$P_1(x),\dots, P_r(x)$ of polynomials of one variable and a collection of numbers
$\mu_1,\dots,\mu_n$, $\mu_1+\dots +\mu_n =0$. That equation is an
analog of the differential equation \Ref{BAE new''}. Under certain
conditions, we show that the union of roots of polynomials
$P_1(x),\dots, P_r(x)$ form a solution to the Bethe ansatz equations
if and only if $P_1(x),\dots, P_r(x)$ and $\mu_1,\dots,\mu_n$
satisfy that polylinear differential equation. 

This fact may be considered as a generalization of Stieltjes' 
 Lemma \ref{slt BAE lemma} (see also  \cite{Sti}, and Sec. 6.8 in \cite{Sz})
 to an arbitrary Kac-Moody algebra.

\subsection{}
In the next paper we plan to extend the results of this paper from $\slth$
to other Lie algebras.

\subsection{}
This work has been started when the fourth author visited
Universit\'e Paul Sabatier in Toulouse in May-June of 2006.
He thanks the university for hospitality.

\section{KZ equations and the Gaudin model}
\label{sec KZ}

\subsection{Lie algebra $\slg_{r+1}$}
\label{sec Lie algebra glg}

Consider the Lie algebra $\glg_{r+1}$ with
standard generators $e_{a,b}$, \ $a,b=1,\dots,r+1$.
Set $h_a=e_{a,a} - e_{a+1,a+1}$ for $a=1,\dots,r$.
Consider the Lie subalgebra
$$
\g = \slg_{r+1}
$$
generated by elements
$e_{a,b}$, $a\neq b$, and $h_a$, $a=1,\dots,r$.
Let
$\g = \n_ - \oplus \h \oplus \n_+$,
\bea
\n_-\ =\ \oplus_{a>b}\ \C\cdot e_{a,b}\ ,
\qquad
\h\ =\ \oplus_{a=1}^{r}\ \C\cdot h_{a}\ ,
\qquad
\n_+\ =\ \oplus_{a<b} \ \C\cdot e_{a,b}\ ,
\eea
be the Cartan decomposition. Let $\al_a \in \h^*,\ a=1,\dots, r$, be simple roots.

Fix the invariant scalar product on $\g$ such that $(h_i,h_i)=2$.
The scalar product identifies $\g$ and its dual $\g^*$.

\subsubsection{}
The Casimir element $\Omega \in \g^{\otimes 2}$ is the element
$\sum_i x_i\otimes x_i$, where $(x_i)$ is an orthonormal basis of $\g$.

\subsubsection{}
For $\g=\slt$,
\bea
\Omega \ =\
e_{2,1}\otimes e_{1,2} +
e_{1,2}\otimes e_{2,1} \ +
\frac 12 h_{1}\otimes h_1 \ .
\eea
For $\g=\slth$,
\begin{align*}
\qquad \Omega \ =\ {}
& e_{2,1}\otimes e_{1,2} +
e_{1,2}\otimes e_{2,1} +{}
\\[4pt]
& e_{3,2} \otimes e_{2,3} +
e_{2,3}\otimes e_{3,2} +
e_{3,1}\otimes e_{1,3} +
e_{1,3}\otimes e_{3,1} +{}
\\[4pt]
& h_{1}\otimes\Bigl(\,\frac 23 h_1 + \frac 13 h_2\Bigr) +
h_{2}\otimes\Bigl(\,\frac 13 h_1 + \frac 23 h_2\Bigr) \ .
\end{align*}

\subsubsection{}
\label{PBW is defined}

Let $\un$ be the universal enveloping algebra of $\n_-$,
\bea
\un\ =\ \oplus_{\bs l \in \Z_{\geq 0}^r}\ \un [\bs l]\ ,
\eea
where for $\bs l = (l_1,\dots,l_r)$,\ the space $\un [\bs l]$
consists of elements $f$ such that
$$
[f,h]\ = \ \langle h\,,\, \sum_{i=1}^r\,l_i\al_i \rangle f\ .
$$
The element $\prod_i e_{a_i,b_i}$ with $a_i>b_i$
belongs to the graded subspace $\un\,[\bs l]$, where
$\bs l = \sum_i \bs l^{(i)}$ with
$$
\bs l^{(i)}\ =\ (0,0, \dots, 0, 1_b,1_{b+1},
\dots, 1_{a-1}, 0,0, \dots, 0)\ .
$$

Fix an order on the set of
elements $e_{a,b}$ with $r+1 \geq a>b\geq 1$.
Set $e_{a,b} < e_{a',b'}$ if $b < b'$ or $b = b'$ and $a < a'$.

For example $e_{2,1} < e_{3,1} < e_{3,2}$.

Then the ordered products
$\prod_{a>b} e_{a,b}^{n_{a,b}}$ form a graded PBW basis of $\un$.

\subsubsection{}
The grading of $\un$ induces the grading of $\un^{\otimes n}$ for any
$n$, \
$$
\un^{\otimes n}\, =\, \oplus_{\bs l}\, \un^{\otimes n}[\bs l]\ ,
$$
where
$\un^{\otimes n}[\bs l] = \oplus_{\bs l^{(1)} + \dots + \bs l^{(n)}= \bs l}
\un[\bs l^{(1)}] \otimes \dots \otimes \un[\bs l^{(n)}].$
The PBW basis of tensor factors induces a graded PBW basis of
$\un^{\otimes n}$.

\subsubsection{}

For $\La \in \h^*$, denote by $M_\La$ the
$\slg_{r+1}$ Verma module with highest weight $\La$. Denote by
$v_\La \in M_\La$ its highest weight vector.

For $\bs \La = (\La_1,\dots,\La_n)$,
$\La_s \in \h^*$, denote
$$
M_{\bs \La}\ =\ M_{\La_1} \otimes \dots \otimes M_{\La_n}\ .
$$
We have the weight decomposition
$$
M_{\bs \La}
\ =\ \oplus_{\bs l\in \Z_{\geq 0}^r}\,M_{\bs \La}[\bs l]\ ,
$$
where $M_{\bs \La}[(l_1, \dots,l_r)]$ is the subspace
of vectors of weight $\sum_{s=1}^n\La_s - \sum_{i=1}^r l_i\al_i$.

The PBW basis of $\un^{\otimes n}$ induces a graded PBW basis of
$M_{\bs \La}$.

Let
$$
\sing M_{\bs \La}[ \bs l]\ \subset \ M_{\bs \La}[\bs l]
$$
denote the subspace
of singular vectors, i.e. the subspace of vectors annihilated by $\n_+$.

\subsection{KZ equations} The KZ equations on an
$M_{\bs \La}$-valued function $I(z_1,\dots,z_n)$ of complex variables
$\bs z = (z_1,\dots,z_n)$ is the following system of
differential equations
\bea
\kappa\ \frac {\partial I}{\partial z_i}\ =\
\sum_{j\neq i}\ \frac{\Omega^{(ij)}}{z_i-z_j}\ I\ ,
\qquad
i = 1, \dots, n\ ,
\eea
where $\kappa \in \C^*$ is a parameter of the equations.

The KZ equations will not be used in this paper, but it is useful to
keep them in mind when discussing the Gaudin model.

\subsection{Gaudin model} Fix pairwise distinct complex numbers
$z_1,\dots,z_n$. Denote
\bea
H_i(\bs z)\ =\ \sum_{j\neq i}\ \frac{\Omega^{(ij)}}{z_i-z_j} \ ,
\qquad
i = 1, \dots, n\ .
\eea
These are linear operators on $M_{\bs \La}$
called {\it the Gaudin Hamiltonians}.

The Gaudin Hamiltonians commute, $[H_i(\bs z), H_j(\bs z)]=0$ for $i\neq j$.

The Gaudin Hamiltonians commute with the action of $\slg_{r+1}$ on
$M_{\bs \La}$ and hence they preserve the subspaces
$\sing M_{\bs \La}[ \bs l]$.

We also have $H_1(\bs z) + \dots + H_n(\bs z) =0$.

\section{Master function, canonical weight function, and the Bethe ansatz}
\label{MasteR}

\subsection{Master function, \cite{SV}}
\label{sec master function}
Let $\bs \La\, =\, (\La_1,\dots,\La_n)$,
$\La_s \in \h^*$, be a collection of
$\slg_{r+1}$-weights and $\bs l = (l_1,\dots,l_r)$
a collection of nonnegative integers.
Set $l = l_1+\dots+l_r$.
Introduce a function of $n$ variables $\bs z=(z_1,\dots,z_n)$ and $l$ variables
\bea
\label{var}
&
\bs t = (t^{(1)}_{1},\dots,t_{l_1}^{(1)},\dots,
t^{(r)}_{1},\dots,t_{l_{r}}^{(r)})\
\eea
by the formula
\bean
\label{master}
&&
\Phi(\bs t;\bs z;\bs \La; \bs l)\ =\ \prod_{1\leq i<j \leq n}
(z_i-z_j)^{(\La_i,\La_j)}\
\times
\\
&&
\phantom{aaa}
\prod_{i=1}^r \prod_{j=1}^{l_i}\prod_{s=1}^n
(t_j^{(i)}-z_s)^{-(\La_s, \al_i)} \prod_{i=1}^r\prod_{1\leq j<s\leq
l_i} (t_j^{(i)}-t_s^{(i)})^{2}
\prod_{i=1}^{r-1}\prod_{j=1}^{l_i}\prod_{k=1}^{l_{i+1}}
(t_j^{(i)}-t_k^{(i+1)})^{-1} \ .
\notag
\eean
The function $\Phi$ is a (multi-valued) function of $\bs t$, depending on
parameters $\bs z, \bs \La$. The function
$\Phi$ is called {\it the master function\/}.

\subsubsection{}
The product of symmetric groups
$$
S_{\bs l}\ =\ S_{l_1}\times \dots \times
S_{l_r}
$$
acts on variables $\bs t$ by permuting the variables with the
same upper index. The master function is $S_{\bs l}$-invariant.

\subsubsection{}
A point $\bs t$ with complex coordinates will be called {\it a critical point}
of $\Phi(\,\cdot\,;\bs z; \bs \La; \bs l)$ if the following system of $l$
equations is satisfied
\bean
\label{critical}
\left(\Phi^{-1}
\frac{\partial \Phi }{\partial t_j^{(i)}}\right)
(\bs t; \bs z; \bs \La; \bs l)\ =\ 0\ ,
\qquad
i = 1 , \dots , r,\ j = 1 , \dots l_i\ .
\eean

\subsubsection{}
Equations \Ref{critical} can be reformulated as the system of equations
\bean\label{BAE}
\sum_{s=1}^{n}
\frac{(\La_s,\al_1)}{t_j^{(1)}-z_s}
-\sum_{s=1,\ s\neq j}^{l_1}
\frac{2}{t_j^{(1)}-t_s^{(1)}}+\sum_{s=1}^{l_2}
\frac{1}{t_j^{(1)}-t_s^{(2)}}=0\ ,
\\
\sum_{s=1}^{n}
\frac{(\La_s,\al_i)}{t_j^{(i)}-z_s}
-\sum_{s=1,\ s\neq j}^{l_i}\frac{2}{t_j^{(i)}-t_{s}^{(i)}}+\sum_{s=1}^{l_{i-1}}
\frac1{t_j^{(i)}-t_s^{(i-1)}}+\sum_{s=1}^{l_{i+1}}
\frac{1}{t_j^{(i)}-t_s^{(i+1)}}=0\ ,
\notag
\\
\sum_{s=1}^{n}
\frac{(\La_s,\al_r)}{t_j^{(r)}-z_s}
-\sum_{s=1,\ s\neq j}^{l_{r}}
\frac{2}{t_j^{(r)}-t_s^{(r)}}+\sum_{s=1}^{l_{r-1}}
\frac{1}{t_j^{(r)}-t_s^{(r-1)}}=0\ ,
\notag
\eean
where $j=1,\dots,l_1$ in the first group of equations,
$i=2,\dots,r-1$ and $j=1,\dots,l_i$ in the second group of
equations, $j=1,\dots,l_{r}$ in the last group of equations.

\subsection{Canonical weight function, \cite{SV, RSV}}

The PBW-basis of $\un$, defined in Section \ref{PBW is defined}, is
$$
\left\{ F_{I_0}\ =\ e_{2,1}^{i_{2,1}} \cdots e_{r+1,r}^{i_{r+1,r}}
\right\}\ ,
$$
where $I_0=\{i_{a,b}\}_{a>b}$ runs over all sequences of nonnegative integers.

Let $F_{I_1},\ldots,F_{I_{n}}$ be elements of the standard
PBW-basis, $I_j=\{i^j_{2,1},\dots, i^j_{r+1,r}\}$. Set
$I=(I_1,\ldots,I_{n})$. The corresponding basis vector
$$
F_Iv\ =\ F_{I_1}v_{\La_1}\otimes\dots\otimes F_{I_{n+1}}v_{\La_n}
$$
lies in $M_{\bs \La}[\bs l]$ if
\begin{equation}\label{pmn}
\sum_{j=1}^{n}\sum_{b=1}^i\sum_{a=i+1}^{r+1} i^j_{a,b}\ =\ l_i\ ,
\qquad
\mbox{ for all }
\quad i=1,\dots,r\ .
\end{equation}
Denote by $P(\bs l,n)$
the set of all indices $I$ corresponding to basis vectors in $M_{\bs \La}[\bs l]$.
The set $\{ F_Iv\}_{I\in P(\bs l,n)}$ forms a basis of $M_{\bs \La}[\bs l]$.

\subsubsection{}
\label{def of canon function}
For $I\in P(\bs l,n)$, define the set
\bea
S(I)\ =\ \{\ (j,a,b,q)\ |\ 1\leq j\leq n,
\quad
1 \leq b < a \leq r+1,
\qquad
1 \leq q \leq i^j_{a,b}\ \}\ .
\eea
For $i=1,\ldots,r$, define the subset
\bea
S_i(I) \ =\ \{ \ s = (j,a,b,q) \in S(I)\ |\ b \leq i <a\ \}\ .
\eea
Condition (\ref{pmn}) implies\ $|S_i(I)|\,=\, l_i$ for $i=1,\ldots,r$.

Define the set
\bea
B(I)\ =\ \{\ \beta = (\beta_1,\dots,\beta_r) \ |\
{\rm for}\ i=1,\dots,r,\
\beta_i\ {\rm is\ a\ bijection} \ S_i(I) \to \{1,\dots,l_i\} \ \}\ .
\eea
We have $|B(I)| = l_1! \cdots l_r!\ .$

For $s=(j,k,l,q) \in S(I)$ and $\beta \in B(I)$, introduce the rational function
\bea
\omega_{s,\beta}\
&=&
\ \frac 1 {t^{(b)}_{\beta_b(s)}-z_j}\
\prod_{i=b+1}^{a-1}\,
\frac 1 {t^{(i)}_{\beta_i(s)}-t^{(i-1)}_{\beta_{i-1}(s)}}\ .
\eea
Introduce the rational functions
\bea
\omega_I\ = \ \frac 1{\prod_{a>b}\ (i^j_{a,b})! }\
\sum_{\beta \in B(I)}\,
\prod_{s\in S(I)}\,
\omega_{s,\beta} \ ,
\qquad
\omega_{\bs l,n}\ =\ \sum_{I \in P(\bs l, n)}\
\omega_I\, F_Iv\ .
\eea
The function $\omega_{\bs l,n}$ defines a rational map
$$
\omega_{\bs l,n}\ :\ \C^l \times \C^n\ \to\ \Ml\ ,
$$
called {\it the canonical weight function}.

\subsubsection{} The canonical weight function was introduced in \cite{SV}.
The formula for the canonical weight function, presented in Section
\ref{def of canon function}, is proved in \cite{RSV}.

\subsubsection{}
{\bf Examples.}
Let $n=2$. If $\bs l = (1, 1, 0, \dots , 0)$, then
\bea
\omega_{\bs l,n}(\bs t,\bs z) &=&
\frac 1{(t_1^{(1)}-z_1)(t_1^{(2)}-z_1)}\ e_{2,1}e_{3,2} v_{\La_1}\otimes v_{\La_2} +
\frac 1{(t_1^{(2)}-t_1^{(1)})(t_1^{(1)}-z_1)}\ e_{3,1}v_{\La_1}\otimes v_{\La_2}
\\
&+&
\frac 1{(t_1^{(1)}-z_1)(t_1^{(2)}-z_2)}\ e_{2,1}v_{\La_1}\otimes e_{3,2}v_{\La_2} +
\frac 1{(t_1^{(2)}-z_1)(t_1^{(1)}-z_2)}\ e_{3,2}v_{\La_1}\otimes e_{2,1}v_{\La_2}
\\
&+&
\frac 1{(t_1^{(1)}-z_2)(t_1^{(2)}-z_2)}\ v_{\La_1}\otimes e_{2,1}e_{3,2}v_{\La_2} +
\frac 1{(t_1^{(2)}-t_1^{(1)})(t_1^{(1)}-z_2)}\ v_{\La_1}\otimes e_{3,1}v_{\La_2}\ .
\eea
If $\bs l = (2, 0, \dots , 0)$, then
\bea
\omega_{\bs l,n}
(\bs t,\bs z)
&=&
\frac 1{(t_1^{(1)}-z_1)(t_2^{(1)}-z_1)}
\ {}e_{2,1}^2v_{\La_1}\otimes v_{\La_2}
\\[3pt]
& +&
\biggl(\frac 1{(t_1^{(1)}-z_1)(t_2^{(1)}-z_2)}+
\frac 1{(t_2^{(1)}-z_1)(t_1^{(1)}-z_2)}\biggr)
\ {} e_{2,1}v_{\La_1}\otimes e_{2,1}v_{\La_2}
\\[2pt]
& +&
\frac 1{(t_1^{(1)}-z_2)(t_2^{(1)}-z_2)}
\ {} \ {}\ v_{\La_1}\otimes e_{2,1}^2v_{\La_2}\ .
\eea

\subsection{Hypergeometric solutions to the KZ equations}

The master function and canonical weight function were introduced in
\cite{SV} to solve the KZ equations.
The hypergeometric solutions to the KZ equations
with values in $\sing \Ml$ have the form
\bea
I(\bs z)\ =\ \int_{\gamma(\bs z)} \ \Phi(\bs t, \bs z, \bs \La, \bs l)^{1/\kappa}
\ \omega_{\bs l, n}(\bs t,\bs z)\ d\bs t\ ,
\eea
where\ $d\bs t\ =\ \wedge_{i,j}\, d t_j^{(i)}$ and $\gamma(\bs z)\ \in\
\C^l\times\{\bs z\}$ is a horizontal family of
$l$-dimensional cycles of
the twisted homology defined by the multi-valued function
$(\Phi)^{\frac1{\kappa}}$, see \cite{SV, V}.

In this paper, we will not use the hypergeometric solutions to the KZ
equations, but the Bethe ansatz for the Gaudin model can be developed
studying quasi-classical asymptotics of these solutions, \cite{RV}.

\subsection{Bethe ansatz for the Gaudin model}
For given $\bs l$ and distinct numbers
$z_1,\dots,z_n$, the problem is to
diagonalize simultaneously the Gaudin Hamiltonians, restricted to
the subspace $\sing \Ml$.

\subsubsection{}
\label{resh}

{\bf Theorem} (\cite{RV}, cf. \cite{Ba1,Ba2}). {\it
Assume that $\bs z \in \C^n$ has distinct coordinates. Assume that
$\bs t \in \C^l$ is a critical point of the master function
$\Phi(\,\cdot\, , \bs z, \bs \La, \bs l)$.
Then the vector $\omega_{\bs l,n}(\bs t,\bs z)$
belongs to $\sing \Ml$ and is an eigenvector of the
Gaudin Hamiltonians with eigenvalues given by the derivatives of the
logarithm of the master function\,:
\bean
\label{eigenvalue}
H_i(\bs z)\
\omega_{\bs l,n}(\bs t, \bs z) \ = \
(\frac{\partial }{\partial z_i}\
{\rm log} \ \Phi (\bs t, \bs z, \bs \La, \bs l))\
\omega_{\bs l,n}(\bs t, \bs z)\ , \qquad i=1,\dots , n\ .
\eean
}

\medskip

This theorem was proved in \cite{RV} using the
quasi-classical asymptotics of the hypergeometric solutions
of the KZ equations. The theorem also follows
directly from Theorem 6.16.2 in \cite{SV},
cf. Theorem 7.2.5 in \cite{SV}.

\subsubsection{}
Notice that a priori the vector $\omega_{\bs l,n}(\bs t, \bs z)$
belongs to $\Ml$, but if $\bs t$ is a critical point, then
$\omega_{\bs l,n}(\bs t, \bs z)$ belongs to $\sing \Ml$.

\subsubsection{}
The critical point equations \Ref{BAE} are called
{\it the Bethe ansatz equations}.

The values of the canonical weight function at the critical points
(with respect to $\bs t$) of the master function are called
{\it the Bethe vectors}.

\section{Different forms of the Bethe ansatz equations for $\slt$ and $\slth$}
\label{sec Different}

\subsection{Bethe ansatz equations for $\slt$}
\label{Bethe ansatz equations for slt}

The Bethe ansatz equations \Ref{BAE} for $\slt$ take the form
\bean\label{BAE2}
\sum_{s=1}^{n}\
\frac{(\La_s,\al_1)}{t_i^{(1)}-z_s}\
-\ \sum_{j=1,\ j\neq i}^{l_1}
\frac{2}{t_i^{(1)}-t_j^{(1)}}\ =\ 0\ ,
\qquad
i=1,\dots,l_1 \ .
\eean

\subsubsection{}
Introduce polynomials
\bea
F(x)\ =\ \prod_{s=1}^n\ (x-z_s)\ ,
\quad
G(x)\ =\ F(x)\, \sum_{s=1}^n\ \frac{(\La_s,\al_1)}{x-z_s}\ ,
\quad
P(x)\! =\! \prod_{i=1}^{l_1}\ (t^{(1)}_i-x)\ .
\eea
We have $\Deg F = n, \ \Deg G = n-1, \ \Deg \,P = l_1$.

\subsubsection{}
\label{slt BAE lemma}
{\bf Lemma.}
{\it
Assume that the roots of $P$ are simple.
Assume that for any $s$ we have
$z_s \notin \{t^{(1)}_1, \dots,t^{(1)}_{l_1}\}$.

Then the roots of $P$ form a solution to
system \Ref{BAE2} if and only if the polynomial
$FP'' - GP'$ is divisible by the polynomial $P$.

In other words,
the roots of $P$ form a solution to the Bethe ansatz equations \Ref{BAE2}
if and only if there exists a polynomial $H$ of degree not greater than
$n-2$ such that $P$ is a solution to the differential equation
\bean
\label{eqn2}
F\,P''\ -\ G\,P'\ +\ H\,P\ =\ 0\ .
\eean
}

The lemma is a  classical result due to Stieltjes,
see \cite{Sti} and Section 6.8 in \cite{Sz}.

\subsubsection{} {\bf Lemma.}
{\it
There exist unique numbers $\mu_1,\dots, \mu_n$,
$\mu_1 +\dots +\mu_n=0$, such that}
\bea
\phantom{aaaaaaaaaaaaaaaaaaaa}
\frac{H(x)}{F(x)}\ =\
\sum_{j=1}^n \frac 1{x-z_j} \,
(\,\mu_j - \sum_{k\neq j}
\frac{ (\La_j,\La_k)}{z_j-z_k}\,)\ .
\phantom{aaaaaaaaaaaaaaaa}
\square
\eea

\subsubsection{}
\label{cor slt}
{\bf Corollary.}
{\it
Assume that the roots of $P$ are simple.
Assume that for any $s$
we have
$z_s \notin \{t^{(1)}_1, \dots,t^{(1)}_{l_1}\}$.
Then the roots of $P$ form a solution to
the Bethe ansatz equations \Ref{BAE2} if and only if there
exist unique numbers $\mu_1,\dots, \mu_n$,
$\mu_1 +\dots +\mu_n=0$, such that
$P$ is a solution to the differential equation}
\bean
\label{eqn2'}
P''\ -\ P' \sum_{i=1}^n\frac{(\La_i,\al_1)}{x-z_i}
\ +\ P \sum_{j=1}^n \frac 1{x-z_j} \,
(\,\mu_j - \sum_{k\neq j}
\frac{ (\La_j,\La_k)}{z_j-z_k}\,)\ = 0\ .
\eean

\subsubsection{}
The corollary provides two ways to describe solutions
to the $\slt$ Bethe ansatz equations.

The original way: a solution is a collection
$\bs t = (t^{(1)}_1, \dots,t^{(1)}_{l_1})$
satisfying \Ref{BAE2}. The second way: a solution is a polynomial $P$
and numbers
$\mu_1,\dots,\mu_n$, $\mu_1+\dots +\mu_n=0$,
satisfying \Ref{eqn2'}.

\subsubsection{}
\label{uniqueness lemma}
{\bf Lemma.}
{\it
For given $F, G, H$, if at least one of the numbers $(\La_1,\al_1),
\dots,$ $(\La_n,\al_1)$ is not an integer, then the polynomial solution
$P$ to equation \Ref{eqn2'} is unique up to multiplication by a nonzero
number.}

Indeed, if $(\La_s,\al_1)$ is not an integer, then $P$ is the unique
(up to multiplication by a number) solution,
uni-valued in a neighborhood of $z_s$.

\subsection{Bethe ansatz equations for $\slth$}

The Bethe ansatz equations \Ref{BAE} for $\slth$ take the form
\bean\label{BAE3}
&&
\sum_{s=1}^{n}
\frac{(\La_s,\al_1)}{t_i^{(1)}-z_s}\
-\ \sum_{j=1,\ j\neq i}^{l_1}
\frac{2}{t_i^{(1)}-t_j^{(1)}}\
+\
\sum_{j=1}^{l_2}
\frac{1}{t_i^{(1)}-t_j^{(2)}}\ =\ 0\ ,
\qquad
i=1,\dots,l_1\ ,
\notag
\\
&&
{}
\\
&&
\sum_{s=1}^{n}
\frac{(\La_s,\al_2)}{t_i^{(2)}-z_s}
\ -\
\sum_{j=1,\ j\neq i}^{l_{2}}
\frac{2}{t_i^{(2)}-t_j^{(2)}}\
+\
\sum_{j=1}^{l_{1}}
\frac{1}{t_i^{(2)}-t_j^{1)}}\ =\ 0\ ,
\qquad
i=1,\dots,l_2\ .
\notag
\eean

\subsubsection{}
Introduce polynomials $F(x) = \prod_{s=1}^n\, (x-z_s)$,\
$P_1(x)\, =\, \prod_{i=1}^{l_1}\, (t^{(1)}_i-x)$,\
\newline
$P_2(x)\, =\, \prod_{i=1}^{l_2}\, (t^{(2)}_i-x)$,\ {}
$F_1(x)\, = \, P_2(x) F(x)$,\ {}
$F_2(x)\, =\, P_1(x) F(x)$,
\bea
G_1(x)\ = \ P'_2(x)F(x)
&+&
P_2(x) F(x) \sum_{s=1}^n \frac{(\La_s,\al_1)}{x-z_s}\ ,
\\
G_2(x)\ =\ P'_1(x)F(x)
&+&
P_1(x) F(x) \sum_{s=1}^n \frac{(\La_s,\al_2)}{x-z_s}\ .
\eea
We have $\Deg\,F_1 = n + l_2, \ \Deg\,F_2 = n + l_1, \
\Deg\,G_1 = n + l_2-1, \ \Deg\,G_2 = n + l_1-1$.

\subsubsection{}
\label{slth BAE lemma}
{\bf Lemma} (\cite{MV1}).
{\it
Assume that the roots
$t^{(1)}_1,\dots,t^{(1)}_{l_1}$,
$t^{(2)}_1,\dots,t^{(2)}_{l_2}$
of $P_1, P_2$ are all pair-wise distinct.
Assume that for any $s, i$,
we have
$z_s \notin \{t^{(i)}_1, \dots,t^{(i)}_{l_i}\}$.

Then the roots of $P_1, P_2$ form a solution to
system \Ref{BAE3} if and only if for $i=1,2$, the polynomial
$F_iP_i'' - G_iP_i'$ is divisible by the polynomial $P_i$.

In other words,
the roots of $P_1, P_2$ form a solution to the Bethe ansatz equations \Ref{BAE3}
if and only if for $i=1,2$,
there exists a polynomial $H_i$ of degree not greater than
$\Deg F_i - 2$ such that $P_i$ is a solution to the differential equation}
\bean
\label{eqn3}
F_i\,P_i''\ -\ G_i\,P_i'\ +\ H_i\,P_i\ =\ 0\ .
\eean

\begin{proof} The lemma is a corollary of the same result
by Stieltjes, see \cite{Sti} and Section 6.8 in \cite{Sz}.
We sketch the proof.

Assume that there exist such polynomials $P_1, P_2, H_1, H_2$.
Substitute $x=t_j^{(i)}$ to equation \Ref{eqn3}. Then we get
$$
\frac{P_i''(t_j^{(i)})}{ P_i'(t_j^{(i)})}\
=\ \frac{G_i(t_j^{(i)})}{ F_i (t_j^{(i)})}\ .
$$
This is exactly the $t^{(i)}_j$-th equation in
\Ref{BAE3}. Hence
the roots of polynomials $P_1,P_2$ form a solution to equations
\Ref{BAE3}.
This argument is reversible.
\end{proof}

\subsubsection{}
\label{new thm}
{\bf Theorem.}
{\it
Assume that the roots
$t^{(1)}_1,\dots,t^{(1)}_{l_1}$,
$t^{(2)}_1,\dots,t^{(2)}_{l_2}$
of $P_1, P_2$ are all pair-wise distinct.
Assume that for any $s, i$,
we have
$z_s \notin \{t^{(i)}_1, \dots,t^{(i)}_{l_i}\}$.
Then the roots of $P_1, P_2$ form a solution to
system \Ref{BAE3} if and only if there exist
numbers $\mu_1, \dots , \mu_n$, $\mu_1 + \dots + \mu_n = 0$,
such that }
\bean
\label{BAE new}
&&
\phantom{aaaaaaa}
P_1''P_2 - P_1'P_2' + P_1P_2'' -
P_1'P_2 \sum_{s=1}^n \frac{(\La_s,\al_1)}{x-z_s}
-
P_1P_2' \sum_{s=1}^n \frac{(\La_s,\al_2)}{x-z_s}
+
\phantom{aaaaaaaaaa}
\\
&&
\phantom{aaaaaaaaaaaaaaaaaaaaaaaaaaaaa}
+ P_1 P_2 \sum_{s=1}^n \frac{1}{x-z_s}
(\,\mu_s - \sum_{k\neq s}
\frac{ (\La_s,\La_k)}{z_s-z_k}\,)\ = 0\ .
\notag
\eean

\subsubsection{}
{\bf Remark.} From \Ref{BAE new} one may conclude that
\bea
\mu_s \ =\ \sum_{k\neq s}
\frac{ (\La_s,\La_k)}{z_s-z_k}\,)
\ +\
(\La_s,\al_1)\frac{P_1'}{P_1}(z_s)
\ +\
(\La_s,\al_2)\frac{P_2'}{P_2}(z_s)
\ =\
\frac{\partial}
{\partial z_s}
{\rm log} \ \Phi (\bs t, \bs z, \bs \La, \bs l))\ ,
\eea
c.f. Theorem \ref{resh}. Thus $\mu_s$ is the eigenvalue of the $s$-th Gaudin operator
at the Bethe vector corresponding to the solution $\bs t$ of the Bethe ansatz equations.

\subsubsection{}
{\it Proof of Theorem \ref{new thm}.}
Let us show that \Ref{BAE new} implies \Ref{BAE3}. Substitute
$x=t^{(1)}_j$ to \Ref{BAE new}. Then
\bea
\frac {P_1''(t^{(1)}_j)}{P_1'(t^{(1)}_j)} -
\frac {P_2'(t^{(1)}_j)}{P_2(t^{(1)}_j)} -
\sum_{s=1}^n
\frac { (\La_s,\al_1) } {t^{(1)}_j-z_s}\ =\ 0\ .
\eea
This is the first of equations in \Ref{BAE3}.
Substituting $x=t^{(2)}_j$ to \Ref{BAE new} we get
the second of equations in \Ref{BAE3}.

Let us show that \Ref{eqn3} imply \Ref{BAE new}. Adding the two equations
of \Ref{eqn3} we get
\bea
&&
F P_1''P_2 - F P_1'P_2'
+
F P_1P_2'' -
F P_1'P_2 \sum_{s=1}^n \frac{(\La_s,\al_1)}{x-z_s}
-
F P_1P_2' \sum_{s=1}^n \frac{(\La_s,\al_2)}{x-z_s}
-
\\
&&
\phantom{aaaaaaaaaaa}
-F P_1'P_2' + H_1P_1 + H_2P_2 = 0\ .
\eea
Equation \Ref{BAE new} will be proved if we show that
$-F P_1'P_2' + H_1P_1 + H_2P_2$ is divisible by $P_1P_2$.
For that it is enough to show that
$-F P_1'P_2' + H_1P_1 + H_2P_2$ is divisible by $P_1$ and
divisible by $P_2$.

From the second of equations in \Ref{eqn3} we get
\bea
-F P_1'P_2' + H_2P_2 =
-FP_1P_2'' + FP_1 P_2'\sum_{s=1}^n \frac{(\La_s,\al_2)}{x-z_s}\ .
\eea
Hence $-F P_1'P_2' + H_1P_1 + H_2P_2$ is divisible by $P_1$.
Similarly it is divisible by $P_2$.
$\square$
\hfill

\subsubsection{}
Lemma \ref{slth BAE lemma} and Theorem \ref{new thm}
provide three ways to describe solutions
to the $\slth$ Bethe ansatz equations.

The original way: a solution is a collection
$\bs t = (t^{(1)}_1, \dots,t^{(1)}_{l_1},t^{(2)}_1, \dots,t^{(2)}_{l_2})$
satisfying \Ref{BAE3}. The second way: a solution is a tuple $P_1$, $P_2$,
$H_1$, $H_2$ satisfying \Ref{eqn3}.
The third way: a solution is a pair $P_1$, $P_2$, and a set of numbers
$\mu_1, \dots , \mu_n$, $\mu_1 + \dots + \mu_n = 0$,
satisfying \Ref{BAE new}.

\subsubsection{} There is a fourth way to describe solutions to
the $\slth$ Bethe ansatz equations. Namely,
under certain conditions, the $S_{\bs l}$-orbits of solutions
are in one to one correspondence with third order linear differential
operators with regular singular points at $z_1,\dots,z_n, \infty$,
with prescribed exponents
at the singular points, and with a quasi-polynomial flag of solutions, see
precise statements in
\cite{MTV}. Under that correspondence to a solution
$\bs t$ one assigns the differential operator
\bea
D_{ \bs t}\ =\ ( \frac{d}{dx} -
\ln' ( \frac { T_{1}T_2 } { P_{2} } ) )\,
( \frac{d}{dx} - \ln' ( \frac {P_{2}T_{1} } {P_{1} } ) )\,
( \frac{d}{dx} - \ln' ( P_1 ) ) \ ,
\notag
\eea
where $\ln'(f)$ denotes $({df}/{dx})/f$ for any $f$, and
$T_i(x) = \prod_{s=1}^n (x-z_s)^{(\La_s,\al_i)}$ for $i=1,2$.

All singular points of $D_{\bs t}$ are regular and lie in
$ \{z_1, \dots , z_n,\infty\}$. The exponents of $D_{\bs t}$ at $z_s$
are
$$
0, (\La_s,\al_1)+1, (\La_s,\al_1+\al_2)+2,
$$
for $s = 1, \dots , n$, and
the exponents of $D_{\bs t}$ at $\infty$ are
$$
-l_1,\ -l_1 -(\sum_{s=1}^n\La_s - l_1\al_1-l_2\al_2, \al_1) - 1,
\
-l_1 -(\sum_{s=1}^n\La_s - l_1\al_1-l_2\al_2, \al_1+\al_2) - 2 .
$$
The differential equation $D_{ \bs t}\,u\,=\,0$ has solutions
$u_1,u_2,u_{3}$ such that
$$
u_1\ =\ P_1\ ,
\qquad
\Wr(u_1,u_2)\ =\ P_2 T_1\ ,
\qquad
\Wr(u_1,u_2, u_3)\ =\ T_1^2T_2\ ,
$$
where $\Wr\,(u_1,\dots, u_i)$ denotes the Wronskian of $u_1,\dots, u_i$.

We will not use this fourth way to describe solutions to the Bethe
ansatz equations in this paper.

\section{On separation of variables for $\slt$}
\label{On separation for slt}

In this section we describe Sklyanin's separation of variables for $\slt$
\cite{Sk}, following exposition in \cite{St}.

\subsection{Polynomial representation} The space of polynomials
$\C[x^{(1)},\dots,x^{(n)}]$ is identified with the tensor product
$$
M_{\bs \La}\ =\ M_{\La_1} \otimes \dots \otimes M_{\La_n}\
$$
of $\slt$ Verma modules by the linear map
\bea
(x^{(1)})^{j^1}\dots (x^{(n)})^{j^n}\ \mapsto\
e_{2,1}^{j^1}v_{\La_1} \otimes \dots \otimes
e_{2,1}^{j^n}v_{\La_n}\ .
\eea
Then the $\slt$ action on $\C[x^{(1)},\dots,x^{(n)}]$
is given by the differential operators,
\bea
e^{(i)} = - x^{(i)}\p_{x^{(i)}}^2 + (\La_i,\al_1) \partial_{x^{(i)}} \ ,
\qquad
h^{(i)} = - 2x^{(i)}\partial_{x^{(i)}} + (\La_i,\al_1)\ ,
\qquad
f^{(i)} = x^{(i)} \ ,
\eea
where $\partial_{x^{(i)}}$ denotes the derivative with respect to
$ x^{(i)}$.

The Gaudin Hamiltonians take the form
\bea
&&
H_i(\bs z) =
\phantom{aaaaaaaaaaaaa}
\\
&&
\phantom{aa}
\sum_{j\neq i}
\frac
{-x^{(i)} x^{(j)}(\p_{x^{(i)}}-\p_{x^{(j)}})^2
+ ((\La_i,\al_1) x^{(j)}-(\La_j,\al_1) x^{(i)})
(\p_{x^{(i)}}-\p_{x^{(j)}}) + (\La_i,\La_j)}
{z_i-z_j}\ ,
\eea
$i = 1, \dots , n$.

\subsection{Change of variables}
\label{change of var slt}
Make the change of variables from
$z_1,\ldots,z_n$, $x^{(1)},\ldots,x^{(n)}$ to
$z_1,\ldots,z_n, u, y^{(1)},\ldots,y^{(n-1)}$
using the relation
\bean
\label{change}
\sum^n_{i=1}\ \frac{x^{(i)}}{t-z_i}\ =\
u\,
\frac{ \prod^{n-1}_{k=1}\, (t - y^{(k)})}
{\prod^n_{i=1}\, (t-z_i)}\ ,
\eean
where $t$ is a variable. This relation defines $u, y^{(1)},
\ldots,y^{(n-1)}$
uniquely up to permutation of $y^{(1)}, \dots, y^{(n-1)}$
unless
$u = \sum_{i=1}^n x^{(i)} =0$.
The map
$$
(z_1,\ldots,z_n, u, y^{(1)},\ldots,y^{(n-1)})\to(z_1,\ldots,z_n,x^{(1)},
\ldots,x^{(n)})
$$
is an unramified covering on the complement to the union of diagonals
$y^{(i)}=y^{(j)}$ and the hyperplane $u=0$.

\subsubsection{}
We have
\bea
x^{(i)}\ &=&\ u\ \frac{\prod_{j=1}^{n-1}\,(z_i - y^{(j)})}
{\prod_{s\neq i}\, (z_i - z_s)} \ ,
\\
\frac{\dpar y^{(j)}}{\dpar x^{(i)}}
\ &=&\
-\ \frac{\prod_{s\ne i}\, (y^{(j)} - z_s)}
{u \, \prod_{l\ne j}
(y^{(j)} - y^{(l)})}
\eea
and then
\bea
\dpar_{x^{(i)}} &=&
\dpar_{u} \ - \ \frac 1 u \ \sum_{j=1}^{n-1} \
\frac{\prod_{s\ne i}\, (y^{(j)} - z_s)}
{ \prod_{l\ne j}
(y^{(j)} - y^{(l)})}\,
\,\dpar_{y^{(j)}} \ ,
\\
\dpar_{x^{(i)}}-\dpar_{x^{(j)}}
&=&
\frac {z_j-z_i} u\
\sum_{l=1}^{n-1} \
\frac{ \prod_{s\notin \{i,j\}}\, (y^{(l)} - z_s)}
{ \prod_{m\ne l}
(y^{(l)} - y^{(m)})}\,
\,\dpar_{y^{(l)}} \ .
\eea

\subsection{Eigenvectors of the Gaudin Hamiltonians}
Assume that we have a common eigenfunction of the Gaudin Hamiltonians
$H_i(\bs z)$ with eigenvalues $\mu_i$.
Then the eigenfunction is annihilated by
the operators $H_i(\bs z) - \mu_i$.

Recall that $H_1(\bs z)+ \dots + H_n(\bs z) = 0$ and hence
$\mu_1+\dots + \mu_n = 0$.

Consider the following operators:
\bea
K_j (\bs z)\ =\ \sum_{i=1}^n\ \frac 1 {y^{(j)} - z_i}\ (H_i(\bs z)-\mu_i)\ ,
\qquad j=1,\dots,n-1\ .
\eea
They annihilate a common eigenfunction of the Gaudin Hamiltonians $H_i(\bs z)$
with eigenvalues $\mu_i$.

\subsubsection{}
\label{Sklyanin thm}
{\bf Theorem} (\cite{Sk}, \cite{St}).
{\it
In variables $u, y^{(1)},\dots,y^{(n-1)}$, we have}
\bea
K_j(\bs z)\ =\ - \p^2_{y^{(j)}}\
+\
\sum_{i=1}^n \, \frac{(\La_i,\al_1)}{y^{(j)}-z_i}\, \p_{y^{(j)}}
\ +\
\sum_{i=1}^n \frac 1{y^{(j)}-z_i} \,
(\,-\mu_i + \sum_{k\neq i}
\frac{(\La_i,\La_k)}{z_i-z_k}\,)\ .
\eea

\medskip
This is the main point of the separation of variables: the operator
$K_j(\bs z)$ depends only on $y^{(j)}$ and does not depend on other
variables $y^{(j')}$ and $u$; this differential operator is the same
operator for all $j$; moreover, it is the same operator as in the
Bethe ansatz equation \Ref{eqn2'}.

\subsubsection{}
More precisely,
Theorem \ref{Sklyanin thm} claims two identities:
\bean
\label{first identity}
\p^2_{y^{(j)}} \ =\
\sum_{i=1}^n \sum_{k\neq i}\
\frac{ x^{(i)} x^{(k)} }
{(y^{(j)} - z_i)(z_i-z_k)}\ (\p_{x^{(i)}}-\p_{x^{(k)}})^2\ ,
\eean
\bean
\label{second identity}
\sum_{i=1}^n \, \frac{(\La_i,\al_1)}{y^{(j)}-z_i}\ \p_{y^{(j)}} \ =\
\sum_{i=1}^n \sum_{k\neq i}\
\frac{(\La_i,\al_1) x^{(k)}-(\La_k,\al_1) x^{(i)}}
{(y^{(j)} - z_i)(z_i-z_k)}\ (\p_{x^{(i)}}-\p_{x^{(k)}})\ .
\eean

\subsection{Canonical weight function} Fix a weight
subspace $M_{\bs \La}[\bs l] \subset M_{\bs \La}$, $\bs l = (l_1)$.
This weight subspace corresponds to
the subspace of $\C[x^{(1)},\dots,x^{(n)}]$ of homogeneous
polynomials of degree $l_1$. The canonical weight function
$\omega_{\bs l, n}$ becomes the following function of
$x^{(1)},\dots,x^{(n)}$, $z_1,\dots,z_n$, $t^{(1)}_1,\dots,t^{(1)}_{l_1}$,
\bea
\prod_{j=1}^{l_1}\,(\,
\sum_{i=1}^n\,
\frac {x^{(i)}} { t^{(1)}_{j}-z_i}\,)\ .
\eea
In variables $u, y^{(1)},\dots, y^{(n-1)}$,
$z_1,\dots,z_n$, $t^{(1)}_1,\dots,t^{(1)}_{l_1}$,
the canonical weight function is
\bean
\label{canon fn sl2}
\Psi(\bs t, \bs z, u, \bs y)\ =\
u^{l_1}\,
\frac{ P(y^{(1)})\dots P(y^{(n-1)})}
{ P(z_1)\dots P(z_{n})}\ ,
\eean
where
$P(x) = \prod_{i=1}^{l_1}(t^{(1)}_i-x)$ as in Section
\ref{Bethe ansatz equations for slt}.

\subsubsection{}
Notice that the canonical weight function, as a function
of $u$, $y^{(1)}$, $\dots,$ $y^{(n-1)}$, is the product
of functions of one variable.
This is another manifestation of separation of variables.

\subsubsection{}
\label{eigenvalue sl2 thm}
{\bf Theorem.}
{\it
Assume that the numbers
$\bs z = (z_1,\dots,z_n)$ are distinct
and the numbers
$\bs t = (t^{(1)}_1, \dots,t^{(1)}_{l_1})$ are distinct.
Assume that for any $s$
we have
$z_s \notin \{t^{(1)}_1, \dots,t^{(1)}_{l_1}\}$.
Assume that for such a  $\bs t$,
the canonical weight function
$\Psi(\bs t, \bs z, u, \bs y)$, as a function of $u, \bs y$,
is an eigenvector of the Gaudin Hamiltonians. Then
$\bs t$ is a solution to the Bethe ansatz equations
\Ref{BAE2}.}

\begin{proof}
If $\Psi(\bs t, \bs z, u, \bs y)$ is an eigenfunction, then it
is annihilated by the operators $K_j(\bs z)$,
$j=1,\dots,n-1$. Hence $\bs t$ is a solution of
\Ref{BAE2} by Corollary \ref{cor slt}.
\end{proof}

\subsubsection{}
Theorem \ref{resh} says that if $\bs t$ is a solution to the Bethe
ansatz equations, then the value at $\bs t$ of the canonical weight
function is an eigenvector of the Gaudin Hamiltonians and is a
singular vector. Theorem \ref{eigenvalue sl2 thm} gives a converse
statement: if the value of the canonical weight function at some point
$\bs t$ is an eigenvector of the Gaudin Hamiltonians, then $\bs t$ is
a solution to the Bethe ansatz equations and that value is a singular
vector.

\subsection{}
\label{lemma on product}
{\bf Lemma.}
{\it
Let at least one of the numbers $(\La_1,\al_1), \dots , (\La_n,\al_1)$
be not an integer. Let $f(y^{(1)},\dots, y^{(n-1)})$ be a polynomial.
Assume that for some $\mu_1,\dots,\mu_n,\
\mu_1+\dots + \mu_n=0$, we have $K_i(z) f = 0$ for $i=1,\dots,n-1$.
Then there exists a polynomial $P$ of one variable such that}
\bea
f(y^{(1)},\dots, y^{(n-1)})\ =\
P(y^{(1)}) \dots P(y^{(n-1)}) \ .
\eea

\begin{proof}
Equation $K_1(\bs z) f = 0$ implies that $f(y^{(1)},\dots,
y^{(n-1)}) = P(y^{(1)})\, g(y^{(2)}, \dots , y^{(n-1)})$, where
$P(y^{(1)})$ is a polynomial and $g(y^{(1)}, \dots , y^{(n-1)})$ a
suitable function. The polynomial $P$ is unique up to multiplication
by a number. Applying the same reasoning to $g$ we get the lemma.
\end{proof}

\section{On separation of variables for $\slth$}
\label{On separation for slth}

\subsection{Polynomial representation}
The space of polynomials
$\C[x_k^{(i)}]^{i=1,\dots,n}_{k=1,2,3}$ of $3n$ variables
is identified with the tensor product
$$
M_{\bs \La}\ =\ M_{\La_1} \otimes \dots \otimes M_{\La_n}\
$$
of $\slth$ Verma modules by the linear map
\bea
&&
(x_1^{(1)})^{j^1_1}\, (x_3^{(1)})^{j^1_3} \,(x_2^{(1)})^{j^1_2}\,
\dots \,
(x_1^{(n)})^{j^n_1}\, (x_3^{(n)})^{j^n_3}\, (x_2^{(n)})^{j^n_2}
\phantom{aaaaaaaaaaaaaaaaaaaaa}
\\
&&
\phantom{aaaaaaaaaaaaaaaaaaaaa}
\mapsto\
e_{2,1}^{j^1_1}\,e_{3,1}^{j^1_3}\,e_{3,2}^{j^1_2}\,
v_{\La_1}
\otimes \ \dots\ \otimes
e_{2,1}^{j^n_1}\,e_{3,1}^{j^n_3}\,e_{3,2}^{j^n_2}\,
v_{\La_n}\ .
\eea
Then the $\slth$ action on $\C[x_j^{(i)}]$
is given by the differential operators,
\bea
e_{2,1}^{(i)} = && x_1^{(i)}\ ,
\qquad e_{3,2}^{(i)} = x_2^{(i)} + x_3^{(i)}\dpar_{x_1^{(i)}}\ ,
\qquad e_{3,1}^{(i)} = x_3^{(i)}\ ,
\\
h_1^{(i)} =& -& 2x_1^{(i)}\p_{x_1^{(i)}} + x_2^{(i)}\p_{x_2^{(i)}} -
x_3^{(i)}\p_{x_3^{(i)}} + (\La_i,\al_1)\ ,
\\
h_2^{(i)} =& -& 2x_2^{(i)}\p_{x_2^{(i)}} +
x_1^{(i)}\p_{x_1^{(i)}} -
x_3^{(i)}\p_{x_3^{(i)}} + (\La_i, \al_2)\ ,
\\
e_{1,2}^{(i)} =& -& x_1^{(i)}\p_{x_1^{(i)}}^2
+ x_{2}^{(i)}\p_{x_1^{(i)}}\p_{x_2^{(i)}} - x_{2}^{(i)}\p_{x_3^{(i)}}
- x_{3}^{(i)}\p_{x_3^{(i)}}\p_{x_1^{(i)}} + (\La_i,\al_1)\p_{x_1^{(i)}}\ ,
\\
e_{2,3}^{(i)}
= &
-& x_{2}^{(i)} \p_{x_2^{(i)}} ^2 + x_{1}^{(i)} \p_{x_3^{(i)}} +
(\La_i,\al_2)\p_{x_2^{(i)}} \ .
\\
e_{1,3}^{(i)} =& -& x_3^{(i)}\p_{x_3^{(i)}}^2 -
x_1^{(i)}\dpar_{x_1^{(i)}}\dpar_{x_3^{(i)}} +
x_2^{(i)}\dpar_{x_1^{(i)}}\dpar_{x_2^{(i)}}^2 -
x_2^{(i)}\dpar_{x_2^{(i)}}\dpar_{x_3^{(i)}} -
\phantom{aaaaaaaaaaaaaaaaa}
\\
&-&
(\La_i,\al_2)\dpar_{x_1^{(i)}}\dpar_{x_2^{(i)}} +
(\La_i,\al_1+\al_2)\dpar_{x_3^{(i)}}\ .
\eea

Then we have the following formula for the Casimir operator:
\bea
&\Omega^{(i,j)}& =
\bigl\{- x^{(i)}_1\dpar_{x_1^{(i)}}^2 + x^{(i)}_2
(\dpar_{x_1^{(i)}}\dpar_{x_2^{(i)}} -
\dpar_{x_3^{(i)}}) - x^{(i)}_3\dpar_{x_1^{(i)}}\dpar_{x_3^{(i)}} +
(\La_i,\al_1)\dpar_{x_1^{(i)}}\bigr\}x_1^{(j)} +
\\
&+& x_1^{(i)}\bigl\{- x^{(j)}_1\dpar_{x_1^{(j)}}^2
+ x^{(j)}_2(\dpar_{x_1^{(j)}}\dpar_{x_2^{(j)}} -
\dpar_{x_3^{(j)}}) - x^{(j)}_3\dpar_{x_1^{(j)}}\dpar_{x_3^{(j)}}
+ (\La_j,\al_1)\dpar_{x_1^{(j)}}\bigr\} +
\\
&+& \bigl\{- x^{(i)}_2\dpar_{x_2^{(i)}}^2 + x^{(i)}_1\dpar_{x_3^{(i)}} +
(\La_i,\al_2)\dpar_{x_2^{(i)}}\bigr\}
\bigl\{x^{(j)}_2 + x^{(j)}_3\dpar_{x_1^{(j)}}\bigr\} +
\\
&+& \bigl\{x^{(i)}_2 + x^{(i)}_3\dpar_{x_1^{(i)}}\bigr\}
\bigl\{- x^{(j)}_2\dpar_{x_2^{(j)}}^2 + x^{(j)}_1\dpar_{x_3^{(j)}}
+ (\La_j,\al_2)\dpar_{x_2^{(j)}}\bigr\} +
\\
&+&
\bigl\{- x^{(i)}_3\dpar_{x_3^{(i)}}^2 -
x^{(i)}_1\dpar_{x_1^{(i)}}\dpar_{x_3^{(i)}}
+ x^{(i)}_2\dpar_{x_1^{(i)}}\dpar_{x_2^{(i)}}^2
- x^{(i)}_2\dpar_{x_2^{(i)}}\dpar_{x_3^{(i)}} -
\phantom{aaaaaaaaaaaaaaaaaaaaaaaa}
\\
&&
\phantom{aaaaaaaaaaaaaaaaaaaaaaaa}
- (\La_i,\al_2)\dpar_{x_1^{(i)}}\dpar_{x_2^{(i)}} +
(\La_i,\al_1+\al_2)\dpar_{x_3^{(i)}}\bigr\}x^{(j)}_3 +
\\
&+& x^{(i)}_3\bigl\{- x^{(j)}_3\dpar_{x_3^{(j)}}^2 -
x^{(j)}_1\dpar_{x_1^{(j)}}\dpar_{x_3^{(j)}}
+ x^{(j)}_2\dpar_{x_1^{(j)}}\dpar_{x_2^{(j)}}^2 -
x^{(j)}_2\dpar_{x_2^{(j)}}\dpar_{x_3^{(j)}} -
\phantom{aaaaaaaaaaaaaaaaaaaaaaaa}
\\
&&
\phantom{aaaaaaaaaaaaaaaaaaaaaaaa}
- (\La_j,\al_2)\dpar_{x_1^{(j)}}\dpar_{x_2^{(j)}}
+ (\La_j,\al_1+\al_2)\dpar_{x_3^{(j)}}\bigr\} +
\eea
\bea
+ \bigl\{- 2x^{(i)}_1\dpar_{x_1^{(i)}} + x^{(i)}_2\dpar_{x_2^{(i)}} -
x^{(i)}_3\dpar_{x_3^{(i)}} + (\La_i,\al_1)\bigr\}
\bigl\{- x^{(j)}_1\dpar_{x_1^{(j)}} - x^{(j)}_3\dpar_{x_3^{(j)}} +
(\La_j,\frac{2\al_1+\al_2}3)\bigr\} +
\eea
\bea
+ \bigl\{- 2x^{(i)}_2\dpar_{x_2^{(i)}} + x^{(i)}_1\dpar_{x_1^{(i)}}
- x^{(i)}_3\dpar_{x_3^{(i)}} + (\La_i,\al_2)\bigr\}
\bigl\{- x^{(j)}_2\dpar_{x_2^{(j)}} - x^{(j)}_3\dpar_{x_3^{(j)}} +
(\La_j,\frac{\al_1+2\al_2}3)\bigr\} \ .
\eea
Rearranging the terms we get
\bean
\label{delta in x variables}
\Omega^{(i,j)}\ = \
\Omega^{(i,j)}_0 + \tilde \Omega^{(i,j)}_0 +
\Omega^{(i,j)}_{>0} + \Omega^{(i,j)}_{<0}\ ,
\eean
where
\bean
\label{delta in x variables 1}
\Omega^{(i,j)}_0 &=& (\La_i,\La_j) +
\bigl(x_1^{(j)}x_2^{(i)}\dpar_{x_2^{(i)}}
- x_1^{(i)}x_2^{(j)}\dpar_{x_2^{(j)}}\bigr)
\bigl(\dpar_{x_1^{(i)}} - \dpar_{x_1^{(j)}}\bigr)\ -
\\
&-&\ x_1^{(i)}x_1^{(j)}\bigl(\dpar_{x_1^{(i)}} - \dpar_{x_1^{(j)}}\bigr)^2
- x_2^{(i)}x_2^{(j)}\bigl(\dpar_{x_2^{(i)}} - \dpar_{x_2^{(j)}}\bigr)^2 \ +
\notag
\\
&+&
\ \bigl((\La_i,\al_1)x_1^{(j)} - (\La_j,\al_1) x_1^{(i)}\bigr)
\bigl(\dpar_{x_1^{(i)}} - \dpar_{x_1^{(i)}}\bigr) \ +
\notag
\\
&+&\
\bigl((\La_i,\al_2)x_2^{(j)} - (\La_j,\al_2)x_2^{(i)}\bigr)
\bigl(\dpar_{x_2^{(i)}} - \dpar_{x_2^{(j)}}\bigr) \ ,
\notag
\eean
\bean
\label{delta in x variables 5}
&&
\phantom{aaaaa}
\tilde \Omega^{(i,j)}_0\ =\
- \bigl(x_3^{(i)}x_1^{(j)} + x_1^{(i)}x_3^{(j)}\bigr)
\bigl(\dpar_{x_1^{(i)}} - \dpar_{x_1^{(j)}}\bigr)
\bigl(\dpar_{x_3^{(i)}} - \dpar_{x_3^{(j)}}\bigr) -
\\
&-&
\ x_3^{(i)}x_3^{(j)}\bigl(\dpar_{x_3^{(i)}} - \dpar_{x_3^{(j)}}\bigr)^2
+ \bigl(x_2^{(j)}x_3^{(i)}\dpar_{x_2^{(j)}}
- x_2^{(i)}x_3^{(j)}\dpar_{x_2^{(i)}}\bigr)
\bigl(\dpar_{x_3^{(i)}} - \dpar_{x_3^{(j)}}\bigr)\ +
\notag
\\
&&
\phantom{aaaaa}
+ \ \bigl((\La_i,\al_1+\al_2)x_3^{(j)}\ -
(\La_j,\al_1+\al_2)x_3^{(i)} \bigr)
\bigl(\dpar_{x_3^{(i)}} - \dpar_{x_3^{(j)}}\bigr) \ ,
\notag
\eean
\bean
\label{delta in x variables 8}
&&
\Omega^{(i,j)}_{>0}\ =\
\bigl(x_2^{(i)}x_3^{(j)}\dpar_{x_2^{(i)}}^2 -
x_3^{(i)}x_2^{(j)}\dpar_{x_2^{(j)}}^2\bigr)
\bigl(\dpar_{x_1^{(i)}} - \dpar_{x_1^{(j)}}\bigr)\ +
\\
&&
\phantom{aaa}
+\
\bigl( (\La_j,\al_2)x_3^{(i)}\dpar_{x_2^{(j)}}
- (\La_i,\al_2)x_3^{(j)}\dpar_{x_2^{(i)}}\bigr)
\bigl(\dpar_{x_1^{(i)}} - \dpar_{x_1^{(j)}}\bigr) \ ,
\notag
\eean
\bean
\label{delta in x variables 10}
\Omega^{(i,j)}_{<0}\ =\
\bigl( x_1^{(i)}x_2^{(j)} - x_1^{(j)}x_2^{(i)}\bigr)
\bigl(\dpar_{x_3^{(i)}} - \dpar_{x_3^{(j)}}\bigr) \ .
\eean
The meaning of this decomposition of the Casimir element will be
explained in Section \ref{Casimir element and degree}.

\subsection{Change of variables}
\label{change of var slth}

Make the change of variables from
$z_1,\ldots,z_n$, $x_1^{(1)}, x_2^{(1)}, x_3^{(1)}$,
$\dots$,
$x_1^{(n)}, x_2^{(n)}, x_3^{(n)}$ to
$z_1,\ldots,z_n$, $u_1,u_2,u_3$, $y_1^{(1)}, y_2^{(1)}, y_3^{(1)},
\ldots, y_1^{(n-1)}, y_2^{(n-1)}, y_3^{(n-1)}$
using the relations
\bean
\label{change slth}
\sum^n_{i=1}\ \frac{x_k^{(i)}}{t-z_i}\ =\
u_k\,
\frac{ \prod^{n-1}_{j=1}\, (t - y_k^{(j)})}
{\prod^n_{i=1}\, (t-z_i)}\ ,
\qquad
k = 1, 2, 3\ ,
\eean
where $t$ is a variable.

Denote $\bs x = (x_1^{(1)}, x_2^{(1)}, x_3^{(1)}$,
$\dots$,
$x_1^{(n)}, x_2^{(n)}, x_3^{(n)})$, \
$\bs u = (u_1,u_2,u_3)$,\
$\bs y = (y_1^{(1)},$ $ y_2^{(1)},$ $ y_3^{(1)},$
\ldots, $y_1^{(n-1)}, y_2^{(n-1)}, y_3^{(n-1)})$.

Relations \Ref{change slth} define $(\bs u,\bs y)$
uniquely up to a permutation of $y^{(j)}_k$'s,
which preserves the lower index,
unless
$u_k = \sum_{i=1}^n x_k^{(i)} =0$ for some of $k$'s.
The map
$$
(\bs z,\bs u,\bs y)\ \to\
(\bs z, \bs x)
$$
is an unramified covering on the complement to the union of diagonals
$y_k^{(i)}=y_k^{(j)}$ and the hyperplanes $u_k=0$.

\goodbreak

\subsection{Degree on $\Ml$}
\label{degree sec}
\subsubsection{}
For an index
$J = (j^1_1,j^1_3,j^1_2 \dots , j^n_1,j^n_3,j^n_2)$ and
$k=1,2,3$,
set $J_k = j^1_k + \dots + j^n_k$.

A monomial
$$
X_J \ =\ (x_1^{(1)})^{j^1_1}\, (x_3^{(1)})^{j^1_3} \,(x_2^{(1)})^{j^1_2}\,
\dots \,
(x_1^{(n)})^{j^n_1}\, (x_3^{(n)})^{j^n_3}\, (x_2^{(n)})^{j^n_2}\
$$
belongs to a weight subspace
$\Ml \subset M_{\bs \La}$, $\bs l = (l_1,l_2)$,
if $J_1 + J_3 = l_1$ and $J_3 + J_2 = l_2$.

We will consider the decomposition
$$
\Ml \ =\ \oplus_{d=0}^{\min\,(l_1,l_2)}
\ M_{\bs \La, d}[\bs l]\ ,
$$
where $M_{\bs \La, d}[\bs l]$ is spanned by all
monomials $X_J$ with
$J_1=l_1-d$, $J_3=d$, $J_2=l_2-d$. We say that
$M_{\bs \La, d}[\bs l]$ consists of {\it
elements of degree $d$.}

\subsubsection{}
In coordinates $(\bs u, \bs y)$, an element belongs to
$M_{\bs \La, d}[\bs l]$, i.e. has degree $d$, if it has the form
\bea
u_1^{l_1-d}\,u_3^{d}\,u_2^{l_2-d}\, f_d\ ,
\eea
where $f_d$ is a polynomial in $\bs y$. This
polynomial does not depend on $y_3^{(1)}$,
$\dots$, $y_3^{(n-1)}$ if $d=0$. An arbitrary
element of $\Ml$ has degree decomposition:
\bea
&&
F\ =\
u_1^{l_1}\,u_2^{l_2}\, f_0(y_1^{(1)},y_2^{(1)},\dots,y_1^{(n-1)},y_2^{(n-1)})
\ +
\phantom{aaaaaaaaaaaaaaaaaaaaaaaaaaaaaaaaaaaaaaaaaaaaaaaaaaaaaaaaaaaaaaaaa}
\\
&&
\phantom{aaaaaaaa}
+ \sum_{d=1}^{\min\,(l_1,l_2)}\,
u_1^{l_1-d}\,u_3^{d}\,u_2^{l_2-d}\,
f_d
(y_1^{(1)},y_3^{(1)},y_2^{(1)},
\dots,y_1^{(n-1)},y_3^{(n-1)},y_2^{(n-1)})\ .
\eea

\subsection{Casimir element and degree}
\label{Casimir element and degree}

The Casimir element $\Omega^{(i,j)}$ is given in
\Ref{delta in x variables}\,--\,\Ref{delta in x variables 10}.
Using formulae of Sections \ref{change of var slth}, consider
the Casimir element as an operator acting on functions of
$(\bs u, \bs y)$.

\subsubsection{}
\label{degree of terms lemma}
{\bf Lemma.}
{\it
\begin{enumerate}
\item[$\bullet$]
The
operator $\Omega^{(i,j)}_0$, given in
\Ref{delta in x variables 1},
preserves the degree introduced in Section \ref{degree sec} and does
not contain
derivatives with respect to
variables $y_3^{(1)}$, \dots, $y_3^{(n-1)}$.
\item[$\bullet$]
The operator $\tilde \Omega^{(i,j)}_0$, given in
\Ref{delta in x variables 5},
preserves the degree and annihilates functions which do not depend
on variables $y_3^{(1)}$, \dots, $y_3^{(n-1)}$.
\item[$\bullet$]
The operator $\Omega^{(i,j)}_{>0}$, given in
\Ref{delta in x variables 8},
increases the degree by one.
\item[$\bullet$]
The operator $\Omega^{(i,j)}_{<0}$, given in
\Ref{delta in x variables 10}, decreases the degree by one
and annihilates functions which do not depend
on variables $y_3^{(1)}$, \dots, $y_3^{(n-1)}$.
\end{enumerate}
}

\medskip
The proof is evident.

\subsubsection{}
In the $3n$-dimensional
space with coordinates
$(\bs u, \bs y)$
consider the subspace
$\frak D$ defined by equations
\bea
y_1^{(j)} \ =\ y_2^{(j)}\ ,
\qquad
j=1,\dots,n-1\ .
\eea
The subspace will be called {\it the main diagonal}.

\subsubsection{}
\label{lem no decreasing}
{\bf Lemma.}
{\it Let $F$ be a polynomial in $(\bs u,\bs y)$. Apply
$\Omega^{(i,j)}_{<0}$ to $F$. 
Then the restriction of the function $\Omega^{(i,j)}_{<0}F$
to the main diagonal equals zero.
}

\medskip

Indeed, the factor $\bigl( x_1^{(i)}x_2^{(j)} -
x_1^{(j)}x_2^{(i)}\bigr)$ in $\Omega^{(i,j)}_{<0}$
is zero on $\frak D$.

\subsubsection{}
\label{cor on the leading term}
Let $F$ be an element in $\Ml$. Let
$F= F_0+F_1+\dots$ and
$\Omega^{(i,j)}F = (\Omega^{(i,j)}F)_0 + (\Omega^{(i,j)}F)_1 + \dots$
be the degree decompositions.

\medskip
\noindent
{\bf Lemma.}
{\it The restrictions to $\frak D$ of the functions
$(\Omega^{(i,j)}F)_0$ and $\Omega^{(i,j)}_0 F_0$
coincide,}
\bea
(\Omega^{(i,j)}F)_0 \vert_{\frak D} \ = \
(\Omega^{(i,j)}_0 F_0 ) \vert_{\frak D} \ .
\eea

\medskip
In other words, the restriction to $\frak D$ of the leading term
$(\Omega^{(i,j)}F)_0$ can be calculated using only the operator
$\Omega^{(i,j)}_0$ applied to $F_0$ and then restricted to $\frak D$.

\begin{proof} The polynomial $F_0$ does not depend on
$y_3^{(1)}, \dots, y_3^{(n-1)}$. Hence
$\tilde \Omega^{(i,j)}_0 F_0 = 0$ and
$\Omega^{(i,j)}_{<0} F_0 = 0$. The function
$\Omega^{(i,j)}_{>0} F_0$ has degree one.
The function
$\Omega^{(i,j)} (F_2 + F_3 + \dots)$ has no degree zero part.
The restriction to $\frak D$ of the degree zero part of the function
$\Omega^{(i,j)} F_1$ is zero by Lemma \ref{lem no decreasing}.
This proves the lemma.
\end{proof}

\subsection{Eigenvectors of the Gaudin Hamiltonians}

\subsubsection{}
Assume that we have a common eigenfunction $F$
of the Gaudin Hamiltonians
$H_i(\bs z)$ with eigenvalues $\mu_i$, $i=1,\dots,n$.
Then the eigenfunction is annihilated by
the operators $H_i(\bs z) - \mu_i$.
Recall that $\mu_1+\dots + \mu_n = 0$.

Consider the following operators:
\bea
K_j (\bs z)\ =\ \sum_{i=1}^n\ \frac 1 {y^{(j)} - z_i}\ (H_i(\bs z)-\mu_i)\ ,
\qquad j=1,\dots,n-1\ .
\eea
They annihilate a common eigenfunction of the Gaudin Hamiltonians $H_i(\bs z)$
with eigenvalues $\mu_i$.

\subsubsection{}
Decomposition \Ref{delta in x variables} of the Casimir operators
into graded components, induces the decomposition of
the Gaudin Hamiltonians into graded components,
\bea
H_i(\bs z)
= \
H_i(\bs z)_0 + \tilde H_i(\bs z)_0 +
H_i(\bs z)_{>0} + H_i(\bs z)_{<0}\ ,
\eea
and the decomposition into graded components
\bea
K_j(\bs z)
= \
K_j(\bs z)_0 + \tilde K_j(\bs z)_0 +
K_j(\bs z)_{>0} + K_j(\bs z)_{<0}\ ,
\eea
where
\begin{alignat*}2
& K_j (\bs z)_0\ =\ \sum_{i=1}^n\, \frac 1
{y_1^{(j)} - z_i}\ (H_i(\bs z)_0-\mu_i)\ ,
\qquad &
\tilde K_j (\bs z)_0\ &{}=\ \sum_{i=1}^n\, \frac 1
{y_1^{(j)} - z_i}\ {\tilde H_i(\bs z)_0} \ ,
\\[4pt]
& K_j (\bs z)_{>0}\ =\ \sum_{i=1}^n\, \frac 1
{y_1^{(j)} - z_i}\ {H_i(\bs z)_{>0}} \ ,
& K_j (\bs z)_{<0}\  &{}=\ \sum_{i=1}^n\, \frac 1
{y_1^{(j)} - z_i}\, {H_i(\bs z)_{<0}} \ .
\end{alignat*}

\subsubsection{}
\label{important equations}
{\bf Lemma.}
{\it
Let $F = F_0+F_1+\dots$ be the degree decomposition of the
common eigenfunction of the Gaudin Hamiltonians, then we have}
\bean
\label{main eqn}
(K_j(\bs z)_0 F_0) \vert_{\frak D}\ =\ 0\ ,
\qquad
j=1,\dots,n-1\ .
\eean

\medskip
The lemma is a corollary of Lemma \ref{cor on the leading term}.

These are important equations. Later, under certain conditions, we
will show how to find $F_0$ from these equations and how to recover
$F$ knowing $F_0$, see Theorem \ref{eigenvalue slth thm}.

\subsubsection{}
\label{thm on separation}
Let $F$ be the
common eigenfunction of the Gaudin Hamiltonians,
$F \in \Ml$, $\bs l =(l_1,l_2)$. Then
\bea
F_0\ =\ u_1^{l_1}\,u_2^{l_2}\, f\ ,
\eea
where $f$ is a polynomial in
$y_1^{(j)}, y_2^{(j)}$, $j=1,\dots,n-1$,
see Section \ref{degree sec}.

\medskip
\noindent
{\bf Theorem.} {\it
Equations \Ref{main eqn} have the form}
\bean
\label{main eqn in coordinates}
&&
\phantom{aa}
\Bigl(\,-\, \frac {\p^2f} {\p y^{(j)}_1 \p y^{(j)}_1}
\ +\
\frac {\p^2f} {\p y^{(j)}_1 \p y^{(j)}_2}
\ -\
\frac {\p^2f} {\p y^{(j)}_2 \p y^{(j)}_2}
\ -\
\frac {\p f} {\p y^{(j)}_1}
\sum_{i=1}^n \frac{(\La_i,\al_1)}{y^{(j)}_1 - z_i} \ -
\phantom{\Bigr)}
\\
&&
\phantom{\Bigl( aaaaaaa}
- \
\frac {\p f} {\p y^{(j)}_2}
\sum_{i=1}^n \frac{(\La_i,\al_2)}{y^{(j)}_1-z_i}
\ +\
f \sum_{i=1}^n \frac 1{y^{(j)}_1-z_i} \,
(\,- \mu_i + \sum_{k\neq i}
\frac{(\La_i,\La_k)}{z_i-z_k}\,)
\,\Bigr) {\bigg{\vert}}_{\frak D} \ =\ 0\ .
\notag
\eean

\medskip

Recall that $\frak D$ is defined by the conditions
$y^{(j)}_1 = y^{(j)}_2$ for $j=1,\dots,n-1$.

\subsubsection{}
As in Sklyanin's Theorem \ref{Sklyanin thm},
equations \Ref{main eqn in coordinates} of Theorem \ref{thm on separation}
have three interesting properties.

The first property is that the variables have separated.
Namely, the $j$-th equation
$(K_j(\bs z)_0 F_0) \vert_{\frak D} = 0$ depends only on variables
$y^{(j)}_1, y^{(j)}_2$ (which at the end are put being equal) and
does not depend
on other variables $y^{(j')}_1$, $y^{(j')}_2$ and $u_1, u_2$.

The second property is that the differential operator
in \Ref{main eqn in coordinates} is
the same equation for all indices $j$.

The third property is that the operator in
\Ref{main eqn in coordinates} is the same as the operator
in the Bethe ansatz
equation \Ref{BAE new} of Theorem \ref{new thm}.

\subsubsection{}
{\it Proof of Theorem \ref{thm on separation}.}
The theorem follows from Theorem \ref{Sklyanin thm}.
Indeed, in order to prove Theorem \ref{thm on separation}
it is enough to prove five identities:
\bean
\label{part1}
\frac {\p^2F_0} {\p y^{(j)}_1 \p y^{(j)}_1}
\ {\bigg{\vert}}_{\frak D}
\ =\
\Bigl(\,
\sum_{i=1}^n\sum_{k\neq i} \
\frac
{ x_1^{(i)} x_1^{(k)} }
{(y^{(j)}_1 - z_i)(z_i-z_k)}
(\dpar_{x_1^{(i)}} - \dpar_{x_1^{(k)}}\bigr)^2
F_0\Bigr)
{\bigg{\vert}}_{\frak D}
\ ,
\notag
\eean
\bean
\label{part2}
\frac {\p^2F_0} {\p y^{(j)}_2 \p y^{(j)}_2}
\ {\bigg{\vert}}_{\frak D}
\ =\ \Bigl(\,
\sum_{i=1}^n\sum_{k\neq i} \
\frac
{ x_2^{(i)} x_2^{(k)} }
{(y^{(j)}_1 - z_i)(z_i-z_k)}
(\dpar_{x_2^{(i)}} - \dpar_{x_2^{(k)}}\bigr)^2
F_0\Bigr)
{\bigg{\vert}}_{\frak D}
\ ,
\notag
\eean
\bean
\label{part3}
&&
\frac {\p^2F_0} {\p y^{(j)}_1 \p y^{(j)}_2}
\ {\bigg{\vert}}_{\frak D}
\ =\
\notag
\\
&&
=\ \Bigl(\,
\sum_{i=1}^n\sum_{k\neq i} \
\frac
{ 1 }
{(y^{(j)}_1 - z_i)(z_i-z_k)}
\bigl(x_1^{(k)}x_2^{(i)}\dpar_{ x_2^{(i)}}
- x_1^{(i)}x_2^{(k)}\dpar_{ x_2^{(k)}}\bigr)
\bigl(\dpar_{x_1^{(i)}} - \dpar_{x_1^{(k)}}\bigr)
F_0\Bigr) {\bigg{\vert}}_{\frak D}
\ ,
\notag
\eean
\bean
\label{part4}
\phantom{aaa}
\frac {\p F_0} {\p y^{(j)}_1}\,
\sum_{i=1}^n \frac{(\La_i,\al_1)}{y^{(j)}_1 - z_i}
\ {\bigg{\vert}}_{\frak D}
\ = \
\Bigl(\,
\sum_{i=1}^n\sum_{k\neq i} \
\frac
{(\La_i,\al_1)x_1^{(k)} - (\La_k,\al_1) x_1^{(i)} }
{(y^{(j)}_1 - z_i)(z_i-z_k)}
(\dpar_{x_1^{(i)}} - \dpar_{x_1^{(k)}}\bigr)
F_0\Bigr)
{\bigg{\vert}}_{\frak D}
\ ,
\notag
\eean
\bean
\label{part5}
\phantom{aaa}
\frac {\p F_0} {\p y^{(j)}_2}\,
\sum_{i=1}^n \frac{(\La_i,\al_2)}{y^{(j)}_1 - z_i}
\ {\bigg{\vert}}_{\frak D}
\ = \
\Bigl(\,
\sum_{i=1}^n\sum_{k\neq i} \
\frac
{(\La_i,\al_2)x_2^{(k)} - (\La_k,\al_2) x_2^{(i)} }
{(y^{(j)}_1 - z_i)(z_i-z_k)}
(\dpar_{x_2^{(i)}} - \dpar_{x_2^{(k)}}\bigr)
F_0\Bigr)
{\bigg{\vert}}_{\frak D}
\ .
\notag
\eean
But the first three identities are corollaries of identity
\Ref{first identity}, and
the last two identities are corollaries of identity \Ref{second identity}.
\hfill
$\square$

\subsection{Canonical weight function} Fix a weight subspace
$\Ml \subset M_{\bs \La}$, $\bs l = (l_1,l_2)$.
In the polynomial representation, the canonical weight function
$\omega_{\bs l,n}$ becomes a function
in $x_1^{(1)},x_2^{(1)},x_3^{(1)}$,
\dots, $x_1^{(n)}, x_2^{(n)}, x_3^{(n)}$,
$z_1,\dots,z_n$, $t^{(1)}_1,\dots,t^{(1)}_{l_1}$,
$t^{(2)}_1,\dots,t^{(2)}_{l_2}$.
After the change of variables of Section \ref{change of var slth}
it becomes a function in $u_1,u_2,u_3$, $y_1^{(1)},y_2^{(1)},y_3^{(1)}$
\dots, $y_1^{(n-1)},
y_2^{(n-1)}, y_3^{(n-1)}$,
$z_1,\dots,z_n$, $t^{(1)}_1,\dots,t^{(1)}_{l_1}$,
$t^{(2)}_1,\dots,t^{(2)}_{l_2}$.

We will give a formula for the canonical weight function and its
degree decomposition. First we prepare notation.

\subsubsection{}
For a given $d=1,\dots, \min\, (l_1,l_2)$, we will sum certain terms
over pairs of ordered subsets $(\bs k, \bs m)$, such that
$\bs k = (k_1,\dots,k_d)$,
$\bs m = (m_1,\dots,m_d)$, where
$k_1,\dots,k_d$ are distinct elements of $(1,\dots,l_1)$
and
$m_1,\dots,m_d$ are distinct elements of $(1,\dots,l_2)$
Such $(\bs k, \bs m)$ will be called {\it $d$-admissible}. The summation over
$d$-admissible pairs will be denoted $\sum_{d-{\rm adm}\,(\bs k, \bs m)}$.

To a $d$-admissible pair $(\bs k, \bs m)$, we assign the following function
\bean
\label{term of can function}
&&
\Xi_{(\bs k, \bs m)}(\bs t, \bs z, \bs y)\ =\
\\
&&
\phantom{aaaa}
=\
\prod_{a=1}^d
\Bigl(
\frac 1{t_{k_a}^{(1)}-t_{m_a}^{(2)}}
\Bigl(
\prod_{j=1}^{n-1}
\frac {t_{k_a}^{(1)}-y_3^{(j)}}
{(t_{k_a}^{(1)}-y_1^{(j)})(t_{m_a}^{(2)}-y_2^{(j)})}
\Bigr)
\prod_{s=1}^n (t_{m_a}^{(2)}-z_s)
\Bigr)
\ .
\notag
\eean

\subsubsection{}
\label{set P's}
Set $P_1(x) = \prod_{i=1}^{l_1} (t^{(1)}_i-x)$,
$P_2(x) = \prod_{i=1}^{l_2} (t^{(2)}_i-x)$.

\subsubsection{}
\label{slth can fn thm}
{\bf Theorem.}
{\it The canonical weight function and its
degree decomposition are given by the following formula,}
\bean
\label{canon fn slth}
&&
\phantom{aa}
\Psi(\bs t, \bs z, \bs u, \bs y)\ =\
u_1^{l_1} u_2^{l_1}\,
\frac{ \prod_{j=1}^{n-1}\,P_1(y_1^{(j)})
P_2(y_2^{(j)})}
{ \prod_{s=1}^n\, P_1(z_s) P_2(z_s)}\,
\times
\\
&&
\phantom{aaaaaaa}
\times
\Bigl[\,
1 + \sum_{d=1}^{\min\,(l_1,l_2)}
\frac 1 {d!}\,
\bigl(\frac {u_3}{u_1 u_2}\bigr)^d \!\!
\sum_{d-{\rm adm}\,(\bs k, \bs m)}\!\!
\Xi_{(\bs k, \bs m)}(\bs t, \bs z, \bs y)\,
\Bigr] \ .
\notag
\eean

\medskip
The theorem is a direct corollary of the definition of the canonical
weight function.

\subsection{Comments on Theorem \ref{slth can fn thm}}
\subsubsection{}
Here are the first terms of the degree decomposition
$$
\Psi(\bs t, \bs z, \bs u, \bs y) = \Psi(\bs t, \bs z, \bs u, \bs y)_0+
\Psi(\bs t, \bs z, \bs u, \bs y)_1 + \Psi(\bs t, \bs z, \bs u, \bs y)_2 + \dots\ .
$$
We have
\bea
\Psi(\bs t, \bs z, \bs u, \bs y)_0\ =\
u_1^{l_1}u_2^{l_1}\
\frac{ \prod_{j=1}^{n-1}\,P_1(y_1^{(j)})
P_2(y_2^{(j)})}
{ \prod_{s=1}^n\, P_1(z_s) P_2(z_s)}\ ,
\eea
\bea
&&
\Psi(\bs t, \bs z, \bs u, \bs y)_1 =
\Psi(\bs t, \bs z, \bs u, \bs y)_0\ \times
\\
&&
\phantom{aaaa}
\times
\frac {u_3}{u_1 u_2}
\sum_{k=1}^{l_1}
\sum_{m=1}^{l_2}
\frac 1{t_{k}^{(1)}-t_{m}^{(2)}}
\Bigl(
\prod_{j=1}^{n-1}
\frac {t_{k}^{(1)}-y_3^{(j)}}
{(t_{k}^{(1)}-y_1^{(j)})(t_{m}^{(2)}-y_2^{(j)})}
\Bigr)
\prod_{s=1}^n (t_{m}^{(2)}-z_s) \ ,
\eea
and so on.

\subsubsection{}
The degree zero term $\Psi_0$ of the canonical weight function
is the product of functions of one variable.
This is a manifestation of separation of variables.

\subsubsection{}
\label{higher terms are determined}
The degree zero term $\Psi_0$ determines all other terms
$\Psi_d$ with $d>0$.

Indeed, knowing $\Psi_0$ we know the roots
$t^{(1)}_1,\dots,t^{(1)}_{l_1}$ of $P_1(x)$ and the roots
$t^{(2)}_1,\dots,t^{(2)}_{l_2}$, of $P_2(x)$.
Now for a $d$-admissible $(\bs k,\bs m)$, the factor
$\Psi_0\, \Xi_{(\bs k,\bs m)}$ has simple combinatorial meaning.
Namely, we cross out $d$ factors with indices in $(\bs k,\bs m)$
from each $P_1, P_2$ entering $\Psi_0$ and then multiply the result by
\bea
\prod_{a=1}^d \
\frac 1{t_{k_a}^{(1)}-t_{m_a}^{(2)}}\
\prod_{j=1}^{n-1} \
\frac{t_{k_a}^{(1)}-y_3^{(j)}}
{\prod_{s=1}^n \,(t_{k_a}^{(1)}-z_s) }
\ .
\eea

Notice that
from $P_1(y^{(j)}_1)$ and $P_2(y^{(j)}_2)$ we cross out factors
depending on $y^{(j)}_1, y^{(j)}_2$ and then multiply the result
by terms
depending on $y^{(j)}_3$. Hence the larger $d$, the more $y^{(j)}_3$-dependent
factors participate in $\Psi_0\, \Xi_{(\bs k,\bs m)}$.

\subsubsection{}
\label{eigenvalue slth thm}
{\bf Theorem.}
{\it
Assume that the numbers $\bs z = (z_1,\dots,z_n)$ are distinct
and the numbers
$\bs t = (t^{(1)}_1, \dots,t^{(1)}_{l_1}$,
$t^{(2)}_1, \dots,t^{(2)}_{l_2})$ are distinct.
Assume that for any $s,i$,
we have
$z_s \notin \{t^{(i)}_1, \dots,t^{(i)}_{l_i}\}$.
Assume that for such a $\bs t$, the canonical weight function
$\Psi(\bs t, \bs z, \bs u, \bs y)$, as a function of $\bs u, \bs y$,
is an eigenvector of the Gaudin Hamiltonians. Then
$\bs t$ is a solution to the Bethe ansatz equations
\Ref{BAE3}.}

\begin{proof}
Let $H_i(\bs z)\Psi(\bs t, \bs z, \bs y, u)= \mu_i
\Psi(\bs t, \bs z, \bs y, u)$ for $i=1,\dots,n$.

Since $\Psi(\bs t, \bs z, \bs y, u)$ is an eigenfunction of the Gaudin
Hamiltonians, it
is annihilated by the corresponding operators $K_j(\bs z)$,
$j=1,\dots,n-1$. Hence the degree zero term
$(K_j(\bs z)\Psi(\bs t, \bs z, \bs y, u))_0$ is zero.
Hence $(K_j(\bs z)_0\Psi(\bs t, \bs z, \bs y, u))_0\vert_{\frak D}$
is zero by Lemma \ref{important equations}.

Using $\bs t$, define $P_1(x), P_2(x)$ as in Section \ref{set P's}.
Then the function $f=P_1(y_1^{(j)}) P_2(y_1^{(j)})$ satisfies equation
\Ref{main eqn in coordinates} by Theorem \ref{thm on separation}.
Then $\bs t$ is a solution to the Bethe ansatz equations \Ref{BAE3}
by Theorem \ref{new thm}.
\end{proof}

\subsubsection{}
Theorem \ref{resh} says that if $\bs t$ is a solution to the Bethe
ansatz equations, then the value at $\bs t$ of the canonical weight
function is an eigenvector of the Gaudin Hamiltonians and is a
singular vector. Theorem \ref{eigenvalue slth thm} gives a converse
statement: if the value of the canonical weight function at some point
$\bs t$ is an eigenvector of the Gaudin Hamiltonians, then $\bs t$ is
a solution to the Bethe ansatz equations and that value is a singular
vector.

\subsection{}
It would be good to have for $\g=\slth$
an analog of Lemma \ref{lemma on product}.

\section{Appendix. New form of the Bethe ansatz
equations for general Kac-Moody algebras}

\subsection{Kac-Moody algebras}\label{Kac_Moody sec}

Let $A=(a_{ij})_{i,j=1}^r$ be a generalized Cartan matrix,
$a_{ii}=2$,
$a_{ij}=0$ if and only $a_{ji}=0$,
$a_{ij}\in \Z_{\leq 0}$ if $i\ne j$.
We assume that $A$ is symmetrizable,
i.e. there is a diagonal matrix $D=\on{diag}\,(d_1,\dots,d_r)$
with positive integers $d_i$ such that $B=DA$
is symmetric.

Let $\g=\g(A)$ be the corresponding complex Kac-Moody
Lie algebra (see \cite{K}, \S 1.2),
$\h \subset \g$ the Cartan subalgebra.
The associated scalar product is non-degenerate on $\h^*$ and
$\on{dim}\h = r + 2d$, where $d$ is the dimension of
the kernel of the Cartan matrix $A$.

Let $\al_i\in \h^*$, $\al_i^\vee\in \h$, $i = 1, \dots , r$, be the
sets of simple roots, coroots, respectively. We have
\bea
(\al_i,\al_j)\ =\ d_i \, a_{ij}\ ,
\qquad
\langle\la ,\al^\vee_i\rangle\
=\
2(\la,\al_i)/{(\al_i,\al_i)}\ {} \ {\rm for}\ {}\
\la\in\h^*\ .
\eea

\subsection{Bethe ansatz equations for the Gaudin model}
Let $\bs z = (z_1,\dots,z_n)$ be a collection of distinct complex numbers.
Let $\bs \La\, =\, (\La_1,\dots,\La_n)$,
$\La_s \in \h^*$, be a collection of
$\g$-weights and $\bs l = (l_1,\dots,l_r)$
a collection of nonnegative integers.
Set $l = l_1+\dots+l_r$ and
\bea
&
\bs t = (t^{(1)}_{1},\dots,t_{l_1}^{(1)},\dots,
t^{(r)}_{1},\dots,t_{l_{r}}^{(r)})\ .
\eea
The Bethe ansatz equations for the Gaudin model associated with this data
is the following system of algebraic equations with respect to $\bs t$:
\bean\label{Bethe eqn}
- \sum_{s=1}^n \frac{(\Lambda_s, \alpha_i)}{t_j^{(i)}-z_s}\ +\
\sum_{s,\ s\neq i}\sum_{k=1}^{l_s} \frac{(\alpha_s, \alpha_i)}{ t_j^{(i)} -t_k^{(s)}}\ +\
\sum_{s,\ s\neq j}\frac {(\alpha_i, \alpha_i)}{ t_j^{(i)} -t_s^{(i)}}\
=\ 0\ ,
\eean
where $i = 1, \dots , r$, $j = 1, \dots , l_i$.

\subsection{The second form of the Bethe ansatz equations}
For a function of $x$ we write $f' = d f/d x$ and $\ln'(f) = f'/f$.

For given $\bs z, \bs \La, \bs t$ and $i=1,\dots,r$, we
introduce polynomials
\bea
P_i(x)\ =\ \prod_{j=1}^{l_i}(t_j^{(i)}-x)\ ,
&&
{T}_i(x)\ =\
\prod_{s=1}^n(x-z_s)^{\langle \La_s, \al_i^\vee\rangle}\ ,
\\
F_i(x)\ = \
\big(\prod_{s=1}^n (x-z_s)\big)\!\!
\prod_{j,\ a_{ij} < 0} P_j(x)\ ,
&&
G_i(x)\ =\ F_i(x) \, \ln'\big(T_i(x)\!\! \prod_{j,\ j\neq i}
P_j(x)^{-\langle\al_j,\al_i^\vee\rangle}\big)
\ .
\eea

\subsubsection{}
\label{deg 2 lem}
{\bf Lemma (\cite{MV1}).}
{\it
Assume that the roots $\bs t$ of polynomials $P_1,\dots,P_r$ are all
simple, distinct and different from $z_1,\dots,z_n$.

Then $\bs t$ is a
solution of the Bethe ansatz equations \Ref{Bethe eqn} if and only if
for every $i = 1, \dots , r$, the polynomial
$F_i P_i'' - G_i P'$ is divisible by the polynomial $P_i$.

In other words, the roots of $P_1,\dots,P_r$ form a solution of the
Bethe ansatz equations \Ref{Bethe eqn} if and only if for every
$i=1,\dots,r$, there exists a polynomial $H_i$ of degree not greater
than $\deg\,F_i - 2$ such that $P_i$ is a solution to the differential equation
}
\bean\label{deg 2 eqn}
F_i\, P_i''\ -\ G_i\, P_i'\ +\ H_i\, P_i\ =\ 0\ .
\eean

\subsection{The new form of the Bethe ansatz equations}
\subsubsection{}
\label{last thm}
{\bf Theorem.}
{\it
Assume that the roots $\bs t$ of polynomials $P_1, \dots, P_r$ are all
simple, distinct and different from $z_1,\dots, z_n$. Then $\bs t$ is a
solution of the Bethe ansatz equations \Ref{Bethe eqn} if and only if
there exist numbers $\mu_1,\dots,\mu_n$,\,
$\mu_1+\dots +\mu_n=0$, such that}
\bean\label{Main}
\sum_{i=1}^r (\al_i,\al_i)\, \frac{P_i''}{P_i} \
+\
\sum_{i\neq j} (\al_i,\al_j) \,\frac{P_i'P_j'}{P_iP_j}-
\sum_{i=1}^r (\al_i,\al_i)\,\frac{T_i'P_i'}{T_iP_i}\ +\
\phantom{aaaaaaaa}
\\
\phantom{aaaaaaaaaaaaaaaaa}
+\
\sum_{s=1}^n\frac{1}{x-z_s}
(\,\mu_j - \sum_{k\neq j}
\frac{ (\La_j,\La_k)}{z_j-z_k}\,)
\ =\ 0\ .
\notag
\eean

This fact may be considered as a generalization of Stieltjes' Lemma
 \ref{slt BAE lemma} (see also  \cite{Sti}, and Sec. 6.8 in \cite{Sz})
 to an arbitrary Kac-Moody algebra.

\begin{proof}
Let
\be
F(x)\ =\
\prod_{s=1}^n (x-z_s)\ ,
\qquad
P(x)\ =\ \prod_{i=1}^rP_i(x)\ .
\ee

Let us show that \Ref{Main} implies \Ref{Bethe eqn}. Multiply
\Ref{Main} by $P$ and substitute for $x$ a root of $P_i$,\,
$x=t_s^{(i)}$. Then
\be
\frac{(\al_i,\al_i)}{2}\frac{P_i''(t_s^{(i)})}{ P_i'(t_s^{(i)})}\
+\
\sum_{j, j\neq i} (\al_i,\al_j) \frac {P_j'(t^{(i)}_s)}{P_j(t^{(i)}_s)}\
-\
\sum_{a=1}^n
\frac { (\La_a,\al_i) } {t^{(i)}_s-z_a}\ =\ 0\ .
\ee
This is exactly the $t^{(i)}_s$-th equation of \Ref{Bethe eqn}.

Let us show the converse. If the zeros of
$P_1,\dots,P_r$ form a solution to \Ref{Bethe eqn}, then for $i=1,\dots,r$ by Lemma
\ref{deg 2 lem}, we have
\bean\label{simple}
\Big(\, \frac{(\al_i,\al_i)}{2} \frac{P_i''}{P_i}
+ \sum_{a, a\neq i}(\al_i,\al_a) \frac{P_i'P_a'}{P_iP_a}-
\frac{(\al_i,\al_i)}{2}\frac{T_i'P_i'}{T_iP_i}\,\Big)FP\ +\ \tilde H_iP_i
\ =\ 0\ ,
\eean
where $\tilde H_i$ is some polynomial.

Add all these equations. Then \Ref{Main} is proved if we show that
\bean
\label{aux}
\sum_{i=1}^r \tilde H_iP_i \ +\ \sum_{i<j}(\al_i,\al_j)
\frac{P_i'P_j'}{P_iP_j} PF
\eean
is divisible by $P$. To show that it suffices to show that
the expression in \Ref{aux} is divisible by $P_1,\dots,P_r$.
Let us show that it is divisible by $P_1$.

Indeed the sum
\bea
\tilde H_1P_1 \ +\ \sum_{1<i<j}(\al_i,\al_j)
\frac{P_i'P_j'}{P_iP_j} PF
\eea
is divisible by $P_1$. Now for $i\neq 1$, from the $i$-th
equation in \Ref{simple} we obtain that
\be
\tilde H_iP_i\ +\ (\al_i,\al_1)
\frac{P_i'P_1'}{P_iP_1} PF
\ee
is divisible by $P_1$. Hence the expression in
\Ref{aux} is divisible by $P_1$.

Note that the quotient of \Ref{aux} by $P$ is a polynomial of degree at most
$n-2$ and therefore the quotient of \Ref{aux} by $PF$ has the form
\bea
\sum_{s=1}^n\frac{1}{x-z_s}
(\,\mu_j - \sum_{k\neq j}
\frac{ (\La_j,\La_k)}{z_j-z_k}\,)
\eea
for some numbers $\mu_1,\dots,\mu_n$ with
$\mu_1+\dots + \mu_n=0$.
\end{proof}

\goodbreak

\end{document}